\documentclass[12pt]{amsart}
\usepackage{amscd,graphicx}

\theoremstyle{definition}

\def \ph{\varphi}

\def \diag{\operatorname {diag}}

\def \ra{\rightarrow}

\def \ie{\hbox{\it i.e.}}

\def \k{\mbox{$\mathbb K$}}

\def \C{\mbox{$\mathbb C$}}
\def \R{\mbox{$\mathbb R$}}
\def \Z{\mbox{$\mathbb Z$}}

\def \P{\mbox{$\mathbb P$}}

\def\zt{\mbox{$\Z_2$}}

\def\inv{^{-1}}
\def\d{d}

\def\im{\operatorname{Im}}
\def\A{\mbox{$\mathcal A$}}

\def\a{\mbox{$\mathfrak a$}}

\def\D{D}

\def\linf{\mbox{$L_\infty$}}
\def\and{\mbox{ \rm and }}

\def\EV{\bigwedge V}
\def\E{\bigwedge}

\def\pha#1#2{\ph^{#1}_{#2}}

\def\psa#1#2{\psi^{#1}_{#2}}

\def\inv{^{-1}}

\def\inf{\text{ inf}}

\author{Alice Fialowski}
\address{E\"otv\"os Lor\'and University\\
Department of Applied Analysis\\
P\'azm\'any P\'eter s\'et\'any 1/C\\
H-1117 Budapest, Hungary}
 \email{fialowsk@cs.elte.hu}
\author{Michael Penkava}
\address{University of Wisconsin\\
Department of Mathematics\\
Eau Claire, WI 54702-4004}
\email{penkavmr@uwec.edu}
\subjclass{14D15,13D10,14B12,16S80,16E40,\\17B56}
\keywords{versal deformations, Lie algebras, moduli space
}
\thanks{The research of the authors was partially supported by
grants from NSF-OTKA-MTA 38453, OTKA T043641 and T043034, and by grants
from the University of Wisconsin-Eau Claire}
\title[Deformations of Four Dimensional Lie Algebras]{Deformations of
Four Dimensional Lie Algebras}

\begin{document}
\setlength{\multlinegap}{0pt}

\begin{abstract}
We study the moduli space of four dimensional
ordinary Lie algebras, and their versal deformations. Their
classification is well
known; our focus in this paper is on the deformations,  which
yield a picture of how the moduli space is assembled. Surprisingly,
we get a nice geometric description of this moduli space essentially 
as an orbifold,
with just a few exceptional points.
\end{abstract}

\maketitle


\section{Introduction}
Lie algebras of small dimension are still a central area of research,
although their classification is basically known up to order 7
(for instance, see \cite{ov,Rh,Tu,KT,pbnl}).
The reason for this is that they play a crucial role in physical
applications (especially in dimension 4). Despite the classification of these
algebras, the moduli space of Lie algebras in a given dimension
is not well understood. We should mention \cite{KN}, on the variety of
$n$ dimensional Lie algebra structures. Moreover, in the existing classifications
there are often overlaps of families determined by parameters and the
manner in which unique objects are singled out is somewhat artificial.
Our solution of this problem is to consider the cohomology of the Lie algebras
as well as their versal deformations, and use this information
as a guide to their division into families. This is the additional
information which provides us with a natural division of the moduli space
of Lie algebras into
families, as well giving us a geometric picture of the structure of the
moduli space.  We did a similar study for 3 dimensional Lie algebras
in \cite{fp3}.

The goal of the present paper is to get an accurate picture of the
moduli space of complex 4-dimensional Lie algebras. The key ingredient in our
description will be the versal deformations of the elements in the moduli space;
therefore cohomology will be a primary computational tool.

In this paper we will
show that the moduli space of Lie algebras on $\C^4$ is essentially an
orbifold given by the natural action of the symmetric group $\Sigma_3$
on the complex projective space $\P^2(\C)$. In addition, there are two
exceptional complex projective lines, one of which has an action of
the symmetric group $\Sigma_2$. Finally, there are 6 exceptional points. The
moduli space is glued together by the miniversal deformations, which
determine the elements that one may deform to locally, so
deformation theory determines the geometry of the space. The
exceptional points play a role in refining the picture of how this
space is glued together.  By orbifold, we mean essentially a topological
space quotiented out by the action of a group. In the case of $\P^n$, there
is a natural action of $\Sigma_{n+1}$  induced by the natural action of
$\Sigma_{n+1}$ on $\C^{n+1}$. An orbifold point is a point which is fixed
by some element in the group. In the case of $\Sigma_{n+1}$ acting on $\P^n$,
points which have two or more coordinates with
the same value are orbifold points, but there are some other ones, such as
the point $(1:-1)=(-1:1)$.

In the classical theory of deformations, a deformation is called a
jump deformation if there is a 1-parameter family of deformations of a
Lie algebra structure such that every nonzero value of the parameter
determines the same deformed Lie algebra, which is not the original
one (see \cite{gers}).
There are also deformations which move along a
family, meaning that the Lie algebra structure is different for each
value of the parameter. There can be multiple parameter families as well.

In the picture we will assemble, both of these phenomena arise.  Some
of the structures belong to families and their deformations simply
move along the family to which they belong.  If there is a jump deformation
from an element
to a member of a family, then there will always be deformations
from that element along the family as well, although they will typically not
be jump deformations.
In addition,
there are sometimes jump deformations either to or from the
exceptional points, so these exceptional points play an interesting
role in the picture of the moduli space.

The structure of this paper is as follows. After some preliminary
definitions and explanation of notation, we will explain our
classification of four dimensional Lie algebras, giving a comparison
between our description of the isomorphism classes of Lie algebras
and the ones in \cite{bs} and \cite{Aga}.  Our division of
the algebras into families is based primarily on cohomological
considerations; elements with the same cohomological description are
placed into the same family in our decomposition. The
correlation between our decomposition and the one in \cite{Aga} is very
close. The main differences arise out of our intention to divide up our
families as projective spaces, a point of view which only partially occurs in \cite{Aga}.

After giving a description of the elements of the moduli space, we
then study in detail miniversal deformations of each element, and
determine how the local deformations behave. The main tool used in this
paper is a constructive approach to the computation of miniversal
deformations, which was first given in \cite{ff2,fp1}. We do not
provide complete details about the method of construction, but try to
provide enough information that the reader might be able to
reconstruct miniversal deformations from the data we provide.  Our
goal here is to use the constructions to give a picture of the moduli
space, rather than to demonstrate the constructions themselves.

Finally, we will assemble all the information we have collected to
give a pictorial representation of the moduli space.

\section{Preliminaries}
In classical Lie algebra theory, the cohomology of a Lie algebra is
studied by considering a differential on the dual space of the
exterior algebra of the underlying vector space, considered as a
cochain complex. If $V$ is the underlying vector space on which the
Lie algebra is defined, then its exterior algebra $\EV$
has a natural
\zt-graded coalgebra structure as well.
In this language, a Lie algebra
is is simply a quadratic odd codifferential on the
exterior coalgebra of a vector space. An odd
codifferential is simply an odd coderivation whose square is zero.
The space $L$ of coderivations has a natural \Z-grading $L=\bigoplus
L_n$, where $L_n$ is the subspace of coderivations determined by
linear maps $\phi:\E^n V\ra V$.  A Lie algebra is a codifferential in
$L_2$, in other words, a quadratic codifferential. (\linf\ algebras are just arbitrary
odd codifferentials.)

The space of
coderivations has a natural structure of a \zt-graded Lie algebra.
The condition that a coderivation $d$ is a codifferential can be
expressed in the form $[d,d]=0$. The coboundary operator $D:L\ra L$ is given
simply by the rule $D(\ph)=[d,\ph]$ for $\ph\in L$; the fact that $D^2=0$ is a direct
consequence of the fact that $d$ is an odd codifferential.  Moreover,
$D(L_n)\subseteq L_{n+1}$, which means that the cohomology $H(d)=\ker
D/\im D $ has a natural decomposition as a \Z-graded space: $H(d)=\prod
H^n(d)$, where
$$
H^n(d)=\ker (D:L_n\ra L_{n+1})/\im (D:L_{n-1}\ra L_n).
$$

Recall that for an arbitrary vector space $V$ of dimension $n$,
the dimension of $\E^k V$ is just $\binom nk$. If $\{e_1,\dots,e_n\}$
is a basis of $V$, and $I=(i_1,\dots,i_k)$ is a multi-index with
$i_1<\cdots<i_k$, and we denote $e_I=e_{i_1}\cdots e_{i_k}$, then the
$e_I$-s give a basis of $\E^k V$. Define $\ph^I_j\in L_k$ by
$\phi^I_j(e_J)=\delta^I_J e_j$, where $\delta^I_J$ is the Kronecker
delta. The elements of $L_k$ are all even if $k$ is odd, and odd if
$k$ is even; to stress this difference, we will denote even
elements as $\phi^I_j$, but odd ones as $\psi^I_j$. Because we will be
working with a four dimensional space, only $L_0$, $L_1$, $L_2$, $L_3$ and
$L_4$ are nonzero, so 1 and 3 cochains are even, while 1, 2 and 4
cochains are odd. In general, the dimension of $L_k$ is just $n\binom
nk$, so for our case, $\dim L_0=4$, $\dim L_1=16$, $\dim L_2=24$, $\dim L_3=16$ and
$\dim L_4=4$.

The Lie algebra structures are codifferentials in $L_2$. In order to
represent a codifferential $d$ as a matrix, we choose the
following order for the increasing pairs $I=(i_1,i_2)$ of indices:
$$\{(1,2),(1,3),(2,3),(1,4),(2,4),(3,4)\},$$ and denote the $i$th element of this ordered set by $S(i)$.
Using this order and the Einstein summation convention, we can
express $$d=a^i_{j}\ph^{S(j)}_i.$$ Let $A=(a^i_j)$ be the matrix
of coefficients of $d$. The first column represents $d(e_1e_2)$,
the second $d(e_1e_3)$, etc.
The Jacobi identity of the Lie algebra is given by the equation $[d,d]=0$, which can be
expressed in matrix form as $AB=0$, where $B$ is the matrix
\begin{equation*}
B :=  \left[
{\begin{array}{cccc}
{a^1_{6}} &  - {a^2_{6}} &  - {a^2_{5}} - {a^1_{4}}
 &  - {a^1_{2}} - {a^2_{3}} \\
 - {a^1_{5}} &  - {a^3_{6}} - {a^1_{4}} &  - {a^3_{5}
} & {a^1_{1}} - {a^3_{3}} \\
 - {a^3_{6}} - {a^2_{5}} &  - {a^2_{4}} & {a^3_{4}}
 & {a^3_{2}} + {a^2_{1}} \\
{a^1_{3}} &  - {a^4_{6}} + {a^1_{2}} &  - {a^4_{5}}
 + {a^1_{1}} &  - {a^4_{3}} \\
 - {a^4_{6}} + {a^2_{3}} & {a^2_{2}} & {a^4_{4}} +
 {a^2_{1}} & {a^4_{2}} \\
{a^4_{5}} + {a^3_{3}} & {a^4_{4}} + {a^3_{2}} &
{a^3_{1}} &  - {a^4_{1}}
\end{array}}
 \right]
\end{equation*}
Since $AB$ is a $4\times 4$ matrix,
we obtain 16 quadratic relations among the coefficients that must be
satisfied. In principle, it should be possible to use a computer
algebra system to determine the solutions,  but in our experience,
this method has some drawbacks, unless one reduces the problem to some
special cases, which we will do below.

In order to classify the solutions,  we note that the dimension of the
derived algebra is just the rank of $A$. We will show that the rank of
$A$ is never larger than 3.  From this it follows that there is an
ideal $I$ of dimension 3 in the Lie algebra $L$, which gives an exact
sequence of Lie algebras
$$
0\ra I\ra L\ra \k\ra 0,
$$
where $\k$ is the abelian Lie algebra of dimension 1. But then, the
structure of $L$ is completely determined by the structure of $I$ as a
Lie algebra, and an outer derivation $\delta$ of $I$.  In \cite{fp3},
the moduli space of three dimensional Lie algebras was studied, and we
will use the classification given there, because we will use in our
classification the
structure of the cohomology of these Lie algebras, which is given in
detail in that paper.

\section{Dimension of the Derived Algebra}
We separate the types of Lie algebras into two distinct cases.
\begin{enumerate}
\item  Every independent pair of vectors spans a two dimensional subalgebra.
\item  There are independent vectors $x$, $y$ and $z$ so that $d(xy)=z$.
\end{enumerate}
The first case is interesting, in that, up to isomorphism, over any field $\k$,
there is
exactly one nonabelian Lie algebra in each dimension greater than one
satisfying this property, and it is given as an extension of a one
dimensional Lie algebra by an abelian ideal. To see this, suppose that
$L$ has dimension at least two, is nonabelian, and satisfies the
property that every independent pair of vectors spans a two
dimensional subalgebra.

Let $x_1'$ and $y'$ be two independent
elements whose bracket $[x_1',y']=ax_1'+by'$ does not vanish. We may
assume that $a\ne 0$.  If $x_1=x_1'+b/ay'$ and $y=1/ay'$, then
$[x_1,y]=x_1$.  Next, suppose that $x_2'$ is independent of $x_1$ and
$y$. Let $[x_2',y]=ax_2'+by$. Then
$x_1+ax_2'+by=[x_1+x_2',y]=p(x_1+x_2')+qy$, for some $p$ and $q$, so
$a=1$. Let $x_2=x_2'+by$. Then $[x_2,y]=x_2$. Now, express $[x_1,x_2]=
ax_1+bx_2$. Then $ax_1+bx_2-x_2=[x_1+y,x_2]=p(x_1+y)+qx_2$, which
implies that $a=0$. Similarly, $x_1+bx_2=[x_1,y+x_2]=px_1+q(y+x_2)$,
so $b=0$ and thus $[x_1,x_2]=0$. The process can be repeated
indefinitely, so we finally obtain a basis $\{x_1,\dots ,x_n,y\}$ satisfying
$[x_i,y]=x_i$, $[x_i,x_j]=0$.

Finally, let us show that the bracket
of any two elements is linearly dependent on them.  Let $u=a^ix_i +by$
and $v=c^ix_i+dy$, then $[u,v]=a^idx_i-c^ibx_i=du-bv$.  Clearly, the
$x_i$-s span an abelian ideal in the algebra, so $L$ is an extension of
the one dimensional Lie algebra spanned by $y$ by this ideal.
It follows that there is an abelian ideal of dimension $n$; moreover, this ideal coincides
with the derived algebra, so the rank of the matrix $A$ is precisely
$n$, one less than the dimension of the vector space.  In fact, the matrix $A$ has precisely
the form $A=\left[\begin{smallmatrix}0&I\\0&0\end{smallmatrix}\right]$, where $I$ is the
$n\times n$ identity matrix. This completes
the description of the first case.

For the second case, suppose that there are linearly independent vectors such that
$d(e_1e_2)=e_3$, so the matrix $A$ of $d$ satisfies
$a^1_{1}=a^2_{1}=a^4_{1}=0$, $a^3_{1}=1$. One can easily check the
possible solutions by considering subcases of this second case. For
example, either $e_1$, $e_2$ and $e_3$ span a subalgebra, or we can
assume that $d(e_1e_2)=e_4$. Since it is well known that the derived
subalgebra of any 4 dimensional Lie algebra has dimension at most 3,
we will not give a detailed analysis of this issue, and simply point out that
the division into subcases can be carried out relatively easily.  However, we note
that even without breaking up the second case into subcases, we can
solve the Jacobi identity using Maple, yielding around 40 solutions
all of which have matrices of rank less than or equal to three. We note that
the solutions are well defined over any field, so the fact that the derived
algebra has dimension 3 is independent of the field $\k$ as well.



\section{Extensions of $\C$ by a three dimensional ideal}
From now on, in this paper, we shall assume that we are working over the base field $\C$.
It is not difficult to classify the moduli space over $\R$ as well. Over fields of
finite characteristic, and over other fields, even the classification of 3 dimensional Lie
algebras is quite complicated.

Since the dimension of the derived algebra is never more than 3, every
4 dimensional Lie algebra is given as an extension of $\C$ by some
three dimensional ideal.  In \cite{fp3}, a complete classification of
three dimensional algebras and their cohomology was given.  We
summarize the results about the cohomology in Table \ref{Table 3}.
\begin{table}[ht]
\begin{center}
\begin{tabular}{llrrrrr}
Type&Codiff&$H^1$&$H^2$&$H^3$\\
\hline\\
$d_1=\mathfrak{n}_3$&$\psi_1^{23}$&4&5&2\\
$d_2=\mathfrak{r}_{3,1}(\C)$&$\psi_1^{13}+\psi_2^{23}$&3&3&0\\
$d_2(1:1)=\mathfrak{r}_{3}(\C)$&$\psi_1^{13}+\psi_1^{23}+\psi_2^{23}$&1&1&0\\
$d_2(\lambda:\mu)=
\mathfrak{r}_{3,\mu/\lambda}(\C)$&$\psi_1^{13}\lambda+\psi_1^{23}+\psi_2^{23}\mu$&1&1&0\\
$d_2(1:0)=\mathfrak{r}_2(\C)\oplus\C$&$\psi_1^{13}+\psi_1^{23}$&2&1&0\\
$d_2(1:-1)=\mathfrak{r}_{3,-1}(\C)$&$\psi_1^{13}+\psi_1^{23}-\psi_2^{23}$&1&2&1\\
$d_3=\mathfrak{sl}_2(\C)$&$\psi_3^{12}+\psi_2^{13}+\psi_1^{23}$&0&0&0\\\\
\hline
\end{tabular}
\end{center}
\caption{\protect Cohomology of Three Dimensional Algebras}\label{Table 3}
\end{table}
Here we have realigned the family of codifferentials as presented in
\cite{fp3} in order to identify elements which have the same
cohomological type as belonging to the same family.  The changes are
actually modest: the family $d_2(\lambda:\mu)$ coincides with
$d(\mu/\lambda)$ of that paper except that the new element $d_2$ was
given as $d(1)$ in the paper, and the element $d(1:1)$ corresponds to
the element $d_2$ in the previous paper. In addition, we have
introduced projective notation for the family $d_2(\lambda:\mu)$.  It
should be noted that $d_2(\lambda:\mu)=d_2(\mu:\lambda)$, so the
family can be identified with $\P^1(\C)/\Sigma_2$, which makes it an
orbifold with orbifold points at $d_2(1:1)$ and $d_2(1:-1)$, where
there is some atypical phenomena in the moduli space. At the point
$d_2(1:1)$, there is a doppelganger $d_2$, whose neighborhoods
coincide with those of the point $d_2(1:1)$, and which also deforms
infinitesimally into $d_2(1:1)$.  At the point $d_2(1:-1)$, there is a
deformation in the $d_3$ direction as well as a deformation in the
direction of the family. Otherwise, members of the family
$d_2(\lambda:\mu)$ deform only in the direction of the family itself.
The codifferential $d_1$ has deformations into every other type of
codifferential except $d_2$, which accounts for why it has such a
large dimension of $H^2$.

In order to determine all the codifferentials of degree 4, it is only
necessary to study the equivalence classes of codifferentials given by
extending $\C$ by a 3 dimensional algebra, via  an outer derivation.
For this reason, in Table \ref{Table 3}, we have denoted by $H^1$ the dimension
of the outer derivations, unlike our convention in \cite{fp3}.
In most cases, an extension of $\C$ by a 3 dimensional algebra is
equivalent to either an extension by the Heisenberg algebra $d_1$, or
an extension by the zero algebra, that is, a three dimensional central
extension of $\C$. For each of the types of 3 dimensional algebras in
our classification in Table \ref{Table 3}, we will analyze the
extensions of $\C$, by studying the outer derivations.

Suppose that $A$ is a matrix representing a codifferential $d$ and
$A'$ is the matrix representing a codifferential $d'$. The
codifferentials $d$ and $d'$ determine isomorphic Lie algebras, and we
call them equivalent codifferentials, if there is a linear
automorphism $g:V\ra V$ such that $d'=g\inv d \tilde g$, where $\tilde
g:\E^2 V\ra \E^2 V$ is the induced isomorphism. If we represent $g$ by
the $4\times4$ matrix $G=(g^i_j)$, where $g(e_j)=g^i_je_i$,
 then $\tilde g$ is represented by
the $6\times6$ matrix $Q$,
in other words, $\tilde g(e_{S(j)})=Q^i_je_{S(i)}$, then
the coefficients of $Q$ are given by the formula
$$
Q^i_j=g^m_kg^n_l-g^m_lg^n_k
, \text{ where } S(i)=(k,l) \text{ and } S(j)=(m,n).$$

It follows that $d$ is equivalent to
$d'$ precisely when there is an invertible matrix $G$ and a corresponding matrix
$Q$ such that $A'=G\inv AQ$.  It is usually easier to check by
computer whether there is a matrix $G$ and corresponding $Q$ so that
$GA'=AQ$, but then one must be careful to check that $\det(G)\ne0$.

\subsection{The simple Lie algebra $d_3=\mathfrak{sl}_2(\C$)}
Since $\mathfrak{sl}_2(\C)$ is simple, all derivations are inner.  As a
consequence, any extension of $\C$ by $\mathfrak{sl}_2(\C)$ is just a direct sum
$\mathfrak{sl}_2(\C)\oplus\C$.  This 4 dimensional algebra is given by the
codifferential
\begin{equation}
d_3=\psi^{12}_3+\psi^{13}_2+\psi^{23}_1,
\end{equation}
which represents the simple algebra $\mathfrak{sl}_2(\C)\oplus\C$ in the BS list
\cite{bs}.
\subsection{The solvable Lie algebra $d_2=r_3(\C)$}
This algebra is given by the codifferential
\begin{equation*}
d_2=\psi^{13}_1+\psi^{23}_2.
\end{equation*}
$H^1(d_2)=\langle\ph^2_2,\ph^1_2,\ph^2_1\rangle$.
Thus a generic outer derivation of $d_2$ is given by
$\delta=\ph^2_2x+\ph^1_2y+\ph^2_1z$.  An extension of $\C$ by $\delta$ is given by
the rule $d(e_ie_4)=\delta(e_i)$. We compute
\begin{equation*}
d(e_1e_4)=e_2y\quad d(e_2e_4)=e_1z+e_2x\quad d(e_3e_4)=0,
\end{equation*}
so that the general formula for an extension $d$ of $\C$ by $d_2$ is
\begin{equation*}
d=\psi^{13}_1+\psi^{23}_2+\psi^{14}_2y+\psi^{24}_1z+\psi^{24}_2x.
\end{equation*}
When $x^2+4yz\ne0$, $d$ is equivalent to
the codifferential
\begin{equation}
d_2^\sharp=\psi^{12}_1+\psi^{34}_3,
\end{equation}
which represents the Lie algebra
$\mathfrak{r}_2\oplus\mathfrak{r}_2$ in the BS list.  When $x^2+4yz=0$ and the three
parameters are not all 0, then the matrix can be transformed into the
matrix of the codifferential
\begin{equation*}
d_1(1:0)=\psa{12}3+\psa{13}3+\psa{23}4,+\psa{14}4,
\end{equation*}
which represents the Lie algebra $\mathfrak g_8(0)$ in the BS list.

\subsection{The solvable algebra $d_2(\lambda:\mu)$}
This algebra is given by the codifferential
$$
d_2(\lambda:\mu)=\psa{13}1\lambda+\psa{23}1+\psa{23}2.
$$
If we consider the trivial extension of $\C$ by $d_2(\lambda:\mu)$,
then $\{e_1,e_2,e_4\}$ span an abelian ideal, so this case reduces to
an extension of $\C$ by an abelian ideal. To analyze nontrivial
extensions, first note that
\begin{align*}
H^1(d_2(1:0))&=\langle\ph^1_1+\ph^2_2,\ph^3_2\rangle\\
H^1(d_2(\lambda:\mu))&=\langle\ph^1_1+\ph^2_2\rangle\quad \text{ otherwise }
\end{align*}
If we extend our codifferential by the derivation
$\delta=(\ph^1_1+\ph^2_2)x+\ph^3_2y$, the extended codifferential is
$$d=d_2(\lambda:\mu)+(\ph^{14}_1+\ph^{24}_2)x+\ph^{34}_2y.$$
When $x\ne 0$ then if  $\lambda=\mu$, the extended codifferential is
equivalent to the codifferential $d_1(1:0)$, otherwise it is
equivalent to the codifferential
$d_2^\sharp$.

When $x=0$ and $\mu\ne0$, the codifferential is equivalent to the
unextended codifferential which we will identify with the
codifferential
\begin{equation*}
d_3(\lambda:\mu:0)=\psa{14}1\lambda+\psa{24}1+\psa{24}2\mu+\psa{34}2,
\end{equation*}
which represents the Lie algebra $\mathfrak
r_{3,\mu/\lambda}(\C)\oplus\C$ (unless $\lambda=\mu$, in which case it
represents the Lie algebra $\mathfrak r_3(\C)\oplus\C$). When $\mu=0$
and $x=0$, then if $y\ne 0$, the extended codifferential is equivalent
to $d_3(1:0:0)$, which represents the Lie algebra $\mathfrak
g_2(0,0)$, but when $y=0$, the unextended codifferential is equivalent
to the codifferential
\begin{equation*}
d_3(0:1)=\psa{34}2+\psa{34}3,
\end{equation*}
which represents the Lie algebra $\mathfrak r_2(\C)\oplus \C^2$.


\subsection{The Heisenberg Algebra $d_1=n_3(\C)$}
Let $d_1=\psi_1^{23}$ be the three dimensional Heisenberg algebra.  Then
\begin{align*}
H^1(d_1)=\langle \ph^2_3,\ph^3_2,\ph^1_1+\ph^2_2,\ph^1_1+\ph^3_3\rangle
\end{align*}
so $H^1(d_1)$ is four dimensional.
If we consider a generic outer derivation
\begin{equation*}
\delta=\ph^2_3a+\ph^3_2b+(\ph^1_1+\ph^2_2)c+(\ph^1_1+\ph^3_3)d,
\end{equation*}
the term
$\psi_3^{24}a+\psi_2^{24}c+\psi_2^{34}b+\psi_3^{34}d+\psi_1^{14}(c+d)$
would be added to $d_1$ obtain the extended codifferential. If we set
$a=a^3_5$, $b=a^2_6$, $c=a^2_5$ and $d=a^3_6$, then we get the
extended codifferential
\begin{equation*}
d=\psi_1^{23}+\psi_1^{14}(a^2_5+a^3_6)+\psi_2^{24}a^2_5+\psi_2^{34}a^2_6+\psi_3^{24}a^3_5+\psi_3^{34}a^3_6,
\end{equation*}
with matrix $A$ given by
$A=
\Bigg[\begin{smallmatrix}
0&0&1&a^2_5+a^3_6&0&0\\
0&0&0&0&a^2_5&a^2_6\\
0&0&0&0&a^3_5&a^3_6\\
0&0&0&0&0&0
\end{smallmatrix}\Bigg]$.
Let $g$ be the linear transformation whose matrix is
$G=\Bigg[\begin{smallmatrix}
1&0&0&0\\
0&p&r&0\\
0&q&s&0\\
0&0&0&1
\end{smallmatrix}\Bigg]$.  Let
$R=\big[\begin{smallmatrix}p&r\\q&s\end{smallmatrix}\big]$, and assume
$\det(R)=1$. Now the matrix $Q$ is given in block form by
$Q=\bigg[\begin{smallmatrix}R&0&0\\0&I&0\\0&0&R\end{smallmatrix}\bigg]$.
The matrix of $d'=g\inv d\tilde g$ is
$A'=
\Bigg[\begin{smallmatrix}
0&0&1&a^2_5+a^3_6&0&0\\
0&0&0&0&a'^2_5&a'^2_6\\
0&0&0&0&a'^3_5&a'^3_6\\
0&0&0&0&0&0
\end{smallmatrix}\Bigg]$,
where
\begin{equation*}
\bigg[\begin{smallmatrix}a'^2_5&a'^2_6\\a'^3_5&a'^3_6\end{smallmatrix}\bigg]=
R\inv \bigg[\begin{smallmatrix}a^2_5&a^2_6\\a^3_5&a^3_6\end{smallmatrix}\bigg]R,
\end{equation*}
which means that if
$V=\bigg[\begin{smallmatrix}a^2_5&a^2_6\\a^3_5&a^3_6\end{smallmatrix}\bigg]$,
then similar submatrices give equivalent codifferentials.  Note that
the $a^1_4$ coefficient $a^2_5+a^3_6$ is just the trace of the matrix
$V$, which is invariant under similarity transformations. Therefore,
looking at the submatrix $V$ alone, we have the following cases
\begin{list}{$\bullet$}
\item
$V=\big[\begin{smallmatrix}\lambda&1\\0&\mu\end{smallmatrix}\big]$,
corresponding to the codifferential
\begin{equation}
d_1(\lambda:\mu)=\psi_1^{23}+\psi_1^{14}(\lambda+\mu)+\psi_2^{24}\lambda+
  \psi_2^{34}+\psi_3^{34}\mu.
\end{equation}
This family of codifferentials should be thought of as a projective
family, parameterizing $\P^1(\C)$. There is an action of $\Sigma_2$ on
this space which identifies $d_1(\lambda:\mu)$ with
$\d_1(\mu:\lambda)$. There are two orbifold points under this action:
$d_1(1:1)$ and $d_1(1,-1)$. We can reasonably expect something unusual
to happen at these orbifold points. In fact, $d_1(1:-1)$ represents
the Lie algebra $\mathfrak{g}_7$ on the BS list while for all other
values, \ie, when $\lambda+\mu\ne 0$, $d_1(\lambda:\mu)$ represents
the Lie algebra
$\mathfrak{g}_8\left(\tfrac{\lambda\mu}{(\lambda+\mu)^2}\right)$.



\item The diagonal matrix $V=\diag(1,1)$. This is the only nonzero
diagonalizable matrix which does not show up in the case above. Its
associated codifferential is given by the formula
\begin{equation}
d_1^\sharp=\psi_1^{23}+2\psi_1^{14}+\psi_2^{24}+\psi_3^{34},
\end{equation}
representing the Lie algebra $\mathfrak{g}_6$.

\item $V=\left[\begin{smallmatrix}0&1\\0&0\end{smallmatrix}\right]$.
Then the extended codifferential is equivalent to
\begin{equation}
d_2^\star=\psa{24}1+\psa{34}2,
\end{equation}
representing the Lie algebra $\mathfrak n_4(\C)$.
\item $V=0$. This is the original, unextended codifferential, which is
equivalent to the codifferential
\begin{equation}
d_1=\psa{24}1,
\end{equation}
representing the Lie algebra $\mathfrak n_3(\C)\oplus\C$.
\end{list}

\subsection{Extensions of $\C$ by an abelian ideal}
Since $H^1(0)=L_1(\C^3)$, the whole 9 dimensional cochain space,
an extension of $\C$ by $\C^3$ is given by a matrix of the form
$
A=
\left[\begin{smallmatrix}
0&0&0&a^1_4&a^1_5&a^1_6\\
0&0&0&a^2_4&a^2_5&a^2_6\\
0&0&0&a^3_4&a^3_5&a^3_6\\
0&0&0&0&0&0
\end{smallmatrix}\right]$.
If we let $V=\Bigg[\begin{smallmatrix}
a^1_4&a^1_5&a^1_6\\
a^2_4&a^2_5&a^2_6\\
a^3_4&a^3_5&a^3_6\\
\end{smallmatrix}\Bigg]$,
then any matrix $V'$ which is similar to $V$ up to multiplication by a
nonzero constant determines an equivalent codifferential.
Since matrices which are constant multiples of each other determine
the same codifferential, we can think of the nonequivalent
codifferentials as being parameterized projectively.  The
decomposition of these matrices into distinct equivalence classes is
as follows.

\begin{list}{$\bullet$}
\item The codifferential
\begin{equation}
d_3(\lambda:\mu:\nu)=\psi_1^{14}\lambda+\psi_1^{24}+\psi_2^{24}\mu+\psi_2^{34}+\psi_3^{34}\nu
\end{equation}
for $(\lambda:\mu:\nu)\in\P^2(\C)/\Sigma_3$, where the action of
$\Sigma_3$ is given by permutation of the coordinates. These points
determine an orbifold with orbifold points occurring along certain
lines ($\P^1(\C)$) where some of the parameters coincide. It might
seem more natural to use diagonal matrices to represent this two
parameter family; the choice here is based on cohomological
considerations.

\item The codifferential
\begin{equation}
d_3(\lambda:\mu)=\psi_1^{14}\lambda+\psi_2^{24}\lambda +\psi_2^{34}+
\psi_3^{34}\mu
\end{equation}
for $(\lambda:\mu)\in\P^1(\C)$. Here there is no action of the symmetric group.
\item The Heisenberg algebra $d_1=\psi_1^{24}$. The only eigenvalue of the matrix is zero,
and it has two Jordan blocks. We will see that every
point in $d_3(\lambda,\mu)$ is infinitesimally close to this point.

\item The solvable algebra $d_2^*$. The matrix has one Jordan block,
with eigenvalue zero.
\item The identity matrix determines the codifferential
\begin{equation}d_3^*=\psi_1^{14}+\psi_2^{25}+\psi_2^{34},\end{equation}
which represents the Lie algebra $\mathfrak g_1(1)$.
\item The zero algebra $d=0$.  Every point is infinitesimally close to
this zero point.
\end{list}
We summarize these results and give the Lie bracket operations in
standard terminology in the table below.
\begin{table}[ht]
\begin{center}
\begin{tabular}{ll}
Type&Brackets\\
\hline
$d_1(\lambda:\mu)$&
$[e_2,e_3]=e_3,[e_1,e_4]=(\lambda+\mu)e_1,$\\&$[e_2,e_4]=\lambda e_2,[e_3,e_4]=e_2+\mu e_3$\\
$d_3(\lambda:\mu:\nu)$&
$[e_1,e_4]=\lambda e_1,[e_2,e_4]=e_1+\mu e_2,[e_3,e_4]=e_2+\nu e_3$\\
$d_3(\lambda:\mu)$&
$[e_1,e_4]=\lambda e_1,[e_2,e_4]=\lambda e_2,[e_3,e_4]=e_2+\mu e_3$\\
$d_1$&$[e_2,e_4]=e_1$\\
$d_1^\sharp$&
$[e_2,e_3]=e_1,[e_1,e_4]=2e_1,[e_2,e_4]=e_2,[e_3,e_4]=e_3$\\
$d_2^*$&$[e_1,e_2]=e_1,[e_3,e_4]=e_2$\\
$d_2^\sharp$&
$[e_1,e_2]=e_1,[e_3,e_4]=e_3$\\
$d_3$&$[e_1,e_2]=e_3,[e_1,e_3]=e_2,[e_2,e_3]=e_1$ \\
$d_3^*$&$[e_1,e_4]=e_1,[e_2,e_4]=e_2,[e_3,e_4]=e_3$\\
\hline
\end{tabular}
\end{center}
\caption{\protect Table of Lie Bracket Operations}\label{Table 5}
\end{table}

\section{Comparison with the Burde-Steinhoff and Agaoka Lists}
The comparison between the Burde-Steinhoff (BS) list and ours
is slightly complicated. On the other hand, our decomposition
is essentially the same as Agaoka's list, so we will just note the corresponding
element, which is of the form ${\mathbf L}_i(\alpha)$ (see \cite{Aga}).

\subsection{$d_1(\lambda:\mu)={\mathbf L}_8(\mu/\lambda)=
\psi_1^{23}+\psi_1^{14}(\lambda+\mu)+\psi_2^{24}\lambda+\psi_2^{34}+\psi_3^{34}\mu$}
\begin{enumerate}
\item When $\lambda+\mu\ne0$, then
$$
d_1(\lambda:\mu)=\mathfrak g_8\left(\tfrac{\lambda\mu}{(\lambda+\mu)^2}\right).
$$
\item When $\lambda+\mu=0$, then we have the codifferential $d_1(1:-1)$ and
$$
d_1(1:-1)=\mathfrak g_7.
$$
\end{enumerate}

\subsection{$d_3(\lambda:\mu:\nu)={\mathbf L}_7(\lambda/\nu,\mu/\nu)
=\psi_1^{14}\lambda+\psi_1^{24}+\psi_2^{24}\mu+\psi_2^{34}+\psi_3^{34}\nu$}
\begin{enumerate}
\item When the trace $\lambda+\mu+\nu$ of the matrix $V$ is nonzero and
none of the parameters are equal to zero, then
$$
d_3(\lambda:\mu:\nu)=
\mathfrak{g}_2\left(\tfrac{\lambda\mu\nu}{(\lambda+\mu+\nu)^3},
\tfrac{\lambda\mu+\lambda\nu+\mu\nu}{(\lambda+\mu+\nu)^2}\right).
$$
\item When exactly one of the parameters vanishes and the other two
are not equal, then
$$
d_3(\lambda:\mu:0)=\mathfrak{r}_{3,\mu/\lambda}(\C)\oplus\C.
$$
\item When one of the parameters vanishes and the other two are equal
we have the special point
$$
d_3(1:1:0)=\mathfrak{r}_3(\C)\oplus\C.
$$
\item When two of the parameters vanish, then we have the special
point
$$
d_3(1:0:0)=\mathfrak{g}_2(0,0).
$$

\item
When the trace of $V$ is zero, none of the parameters is equal to
zero, and the parameters are not the three distinct roots of unity,
then we have
$$
d_3(\lambda:\mu:-\lambda-\mu)=
\mathfrak{g}_3\left(\tfrac{(\lambda^2+\lambda\mu+\mu^2)^3}{(\lambda\mu(\lambda+\mu))^2}\right).
$$
\item When $\lambda$, $\mu$ and $\nu$ are the three distinct cube
roots of unity, then
$$
d_3(1:-1/2+1/2i\sqrt{3}:-1/2-1/2i\sqrt{3})=\mathfrak{g}_4.
$$
\end{enumerate}

\subsection{$d_3(\lambda:\mu)=
{\mathbf L}_4(\mu/\lambda)=
\psi_1^{14}\lambda+\psi_2^{24}\lambda+\psi_2^{34}+\psi_3^{34}\mu$}
\begin{enumerate}
\item When neither of the parameters  vanish or are equal, then we have
$$
d_3(\lambda:\mu)=\mathfrak{g}_1(\mu/\lambda).
$$

\item When $\mu=0$, then we have the special point
$$
d_3(1:0)=\mathfrak{r}_{3,1}(\C)\oplus\C.
$$
\item When $\lambda=0$ then we have the special point
$$
d_3(0:1)={\mathbf L}_4(\infty)=\mathfrak{r}_2(\C)\oplus\C^2.
$$
\item When $\lambda=\mu$ then we have the special point
$$
d(0)=d_3(1:1)=\mathfrak{g}_5.
$$
\end{enumerate}

\subsection{The special cases}
\begin{align*}
d_1&={\mathbf L}_1=\mathfrak{n}_3(\C)\oplus\C=\psi_1^{24}\\
d_1^\sharp&={\mathbf L}_5=\mathfrak g_6=\psi_1^{23}+2\psi_1^{14}+\psi_2^{24}+\psi_3^{34}\\
d_2^*&={\mathbf L}_2=\mathfrak{n}_4(\C)=\psi_1^{24}+\psi_2^{34}\\
d_2^\sharp&={\mathbf L}_9=\mathfrak r_2(\C)\oplus\mathfrak r_2(\C)=\psi_1^{12}+\psi_3^{34}\\
d_3&={\mathbf L}_6=\mathfrak{sl}_2(\C)\oplus\C=\psi_3^{12}+\psi_2^{13}+\psi_1^{23}\\
d_3^*&={\mathbf L}_3=\mathfrak{g}_1(1)=\psi_1^{14}+\psi_2^{24}+\psi_3^{34}.
\end{align*}

\section{Deformations of the Lie Algebras}
For the basic notion of deformations, we refer to \cite{gers,nr,fi1,fi2,ff2}.
In some previous papers, we considered deformations of \linf\ algebras
\cite{fp1, fp3}.
In this paper,
we will only consider Lie algebra
deformations of our Lie algebras, which are determined
by cocycles coming from $H^2$.
We will not explore \linf\
deformations of the Lie algebras we study in this paper, but it would not be difficult to construct
them from the cohomology computations we provide
here.

In Table \ref{Table 2}, we give a classification of the
codifferentials according to their cohomology. Note that for the most
part, elements from the same family have the same cohomology. In fact,
the decomposition of the codifferentials into families was strongly
influenced by the desire to associate elements with the same pattern
of cohomology in the same family.  This is why our family
$d_3(\lambda:\mu:\nu)$ was not chosen to be the diagonal matrices.
Similar considerations influenced our selection of the family
$d_3(\lambda:\mu)$.

\begin{table}[ht]
\begin{center}
\begin{tabular}{lrrrrr}
Type&$H^1$&$H^2$&$H^3$&$H^4$\\
\hline
$d_3$&1&0&1&1\\
$d_2^\sharp$&0&0&0&0\\
\hline
$d_1(1:-1)$&2&2&2&1\\
$d_1(1:0)$&1&2&1&0\\
$d_1(\lambda:\mu)$&1&1&0&0\\
$d_1^\sharp$&3&3&0&0\\
\hline
$d_3(1:-1:0)$&3&5&5&2\\
$d_3(\lambda:\mu:\lambda+\mu)$&2&3&1&0\\
$d_3(\lambda:\mu:0)$&3&3&1&0\\
$d_3(\lambda:\mu:-\lambda-\mu)$&2&2&1&1\\
$d_3(\lambda:\mu:\nu)$&2&2&0&0\\
\hline
$d_3(1:0)$&5&7&3&0\\
$d_3(0:1)$&6&6&2&0\\
$d_3(1:2)$&4&5&1&0\\
$d_3(1:-2)$&4&4&1&1\\
$d_3(\lambda:\mu)$&4&4&0&0\\
\hline
$d_1$&8&13&10&3\\
$d_2^*$&4&6&5&2\\
$d_3^*$&8&8&0&0\\
\hline
\end{tabular}
\end{center}
\caption{\protect Table of the Cohomology}\label{Table 2}
\end{table}

\subsection{The codifferential $d_3=\mathfrak{sl}_2(\C)\oplus\C$}
It is an easy calculation to show that
\begin{align*}
H^1&=\langle\ph^4_4\rangle\\
H^2&=0\\
H^3&=\langle\ph^{123}_4\rangle\\
H^4&=\langle\ph^{1234}_4\rangle
\end{align*}
Since $H^2$ vanishes, this algebra is rigid in terms of deformations
in the Lie algebra sense.
%

\subsection{The codifferential
$d_2^\sharp=\mathfrak{r}_2(\C)\oplus\mathfrak{r}_2(\C)$}
Since the cohomology vanishes entirely, this algebra has
no interesting deformations or extensions.
This algebra is the only four
dimensional Lie algebra which is truly rigid in the $\linf$ algebra
sense, although $d_3$ is also rigid in the Lie algebra sense.

\subsection{The codifferential $d_1(\lambda:\mu)$}
In the generic case we have
\begin{align*}
H^1&=\langle
2\pha11+\pha22+\pha33\rangle\\
H^2&=\langle\psa{14}1+\psa{24}2\rangle
\end{align*}
and all higher cohomology vanishes. Thus, generically, the
infinitesimal deformation is given by
\begin{equation}
d^\inf=d_1(\lambda+t,\mu).
\end{equation}
Since $d^\inf$ is actually a member of the family $d_1(\lambda:\mu)$, it is
immediate that $[d^\inf,d^\inf]=0$, so the infinitesimal deformation is the miniversal
deformation $\d^\infty$ and
the base of the miniversal deformation is $\C[[t]]$. Moreover, it is transparent in
this case that the deformations run in the direction of the family.

\subsection{The codifferential $d_1(1:-1)$}
For this special case there are more
cohomology classes than in the generic case.  We
have
\begin{align*}
H^1&=\langle
2\pha11+\pha22+\pha33,\pha41\rangle\\
H^2&=\langle\psi_1=\psa{14}1+\psa{24}2,\psi_2=\psa{23}4\rangle\\
H^3&=\langle\pha{123}4,\pha{123}1+\pha{234}4\rangle\\
H^4&=\langle\psa{1234}4\rangle.
\end{align*}
Consider the universal infinitesimal deformation
\begin{equation*}
d^\inf=d_1(1:-1)+\psi_1t^1 +\psi_2 t^2.
\end{equation*}
Then we have $\tfrac12[d^\inf,d^\inf]=-(\pha{123}1+\pha{234}4)t^1t^2$, which is a nontrivial cocycle.
It follows that the infinitesimal deformation is miniversal,
and the base of the miniversal deformation is $\A=\C[[t^1,t^2,t^3]]/(t^1t^2)$.

When $t^1t^2\ne0$, the miniversal deformation $d^\infty=d^\inf$ does not
correspond to an actual deformation. The cohomology class of the cocycle
$(\pha{123}1+\pha{234}4)$ is called an obstruction to the extension of the infinitesimal
deformation to higher order.  In order to obtain an actual deformation out of the
miniversal deformation, we need to restrict ourselves to the lines
$t^1=0$ or $t^2=0$, along which the obstruction term vanishes. The cohomology $H^2$,
which gives the tangent space to the moduli space, has
dimension 2, but the deformations actually lie on two curves.
Thus the dimension of the tangent
space does not reveal the complete situation in terms of the deformations;
one needs to construct the versal deformation to get the true picture.

A deformation along the line $t^2=0$ gives $d_1(1+t_1,:-1+t_1)$, the
same pattern as we observed generically.
On the other hand, a deformation along the line $t^1=0$ yields a coderivation
which is equivalent to the codifferential
$d_3$
.
This is an example of a jump deformation,
because $d_1(1:-1)+\psi_2t^2\sim d_3$ for all values of $t^2$. In the classical language
of Lie brackets, we get the following bracket table:
\begin{equation*}
[e_1,e_3]=e_1+t^2e_4,\quad[e_2,e_4]=e_2,\quad[e_3,e_4]=e_2-e_3
\end{equation*}

Both orbifold points in the family $d_1(\lambda:\mu)$ have some
unusual features. The point $d_1(1:-1)$, because it has a jump
deformation out of the family to $d_3$, and the point $d_1(1:1)$
because, as we will see shortly, there is a jump deformation to it
from the element $d_1^\sharp$, which lies outside of the family. What is
surprising is that the point $d_1(1:0)$, which is not an orbifold point,
is also special.

\subsection{The codifferential $d_1(1:0)$}

The cohomology $H^1$ is the same as the generic case, while $H^2$ and $H^3$ are given by
\begin{align*}
H^2&=\langle\psi_1=\psa{14}1+\psa{34}3,\psi_2=\psa{13}2\rangle\\
H^3&=\langle\phi=\pha{134}2\rangle.
\end{align*}
The universal infinitesimal deformation $d^\inf=d_1(1:0)+\psi_i t^i$
satisfies $\tfrac12[d^\inf,d^\inf]=-2\phi t^1t^2$, so it is miniversal and the base
of the versal deformation is $\C[[t^1,t^2]]/(t^1t^2)$.  Along
the line $t^2=0$, $d^\infty=d_1(1,t^1)$, so we deform along the family as in the
generic case.

Along the line $t^1=0$, the deformation $d^\inf$ is equivalent to $d_2^\sharp$ for
all values of $t^2$.
Thus $d_1(1:0)$ has a jump deformation to the element $d_2^\sharp$.  The classical form
of the Lie brackets for the case $t^1=0$ is
\begin{equation*}
[e_1,e_3]=t^2e_2, [e_2,e_3]=e_1,[e_1,e_4]=e_1,[e_2,e_4]=e_2,[e_3,e_4]=e_2
\end{equation*}

\subsection{The codifferential $d_1^\sharp$}
The cohomology is given by
\begin{align*}
H^1&=\langle\pha11+\pha22,\pha23,\pha32\rangle\\
H^2&=\langle\psi_1=\psa{24}3,\psi_2=\psa{34}2,\psi_3=\psa{14}1+\psa{34}3\rangle
\end{align*}
The universal infinitesimal deformation $d^\inf=d_1^\sharp+\psi_it^i$ is
miniversal  since $[d^\inf,d^\inf]=0$, so the
base of the miniversal deformation is just $\A=\C[[t^1,t^2,t^3]]$.

Now let us consider which codifferential $d^\infty=d^\inf$ is equivalent
to.  Even though the deformation is defined for all values of the parameters,
which element we deform to depends in a complicated manner on the parameters.

Except on the plane $t^1=0$, we have $d^\infty\sim d_1(\alpha:\beta)$ where
$$(\alpha,\beta)=2+t^3\pm\sqrt{(t^3)^2+4t^1t^2}.$$

On the plane $t^1=0$, we have $d^\infty\sim d_1(1+t^3:1)$.
In particular, if $t^3=0$, we have $d^\infty\sim d(1:1)$. In fact,
along the entire surface given by $(t^3)^2+4t^1t^2=0$, we have $d^\infty\sim
d(1:1)$, so there is a two parameter family of jump deformations to $d_1(1:1)$.

Thus $d_1^\sharp$ has a jump deformation to $d_1(1:1)$ and
deforms along the family $d_1(\alpha:\beta)$ as if it
were the element $d_1(1:1)$ in this family. This is a pattern which will always
emerge: \emph{If a codifferential has a jump deformation to another codifferential,
then it will deform also to every codifferential to which the element it jumps to
deforms.}

We give the classical form
of the Lie brackets for $d^\infty$:
\begin{equation*}
[e_2,e_3]\!=\!e_1,[e_1,e_4]\!=\!(2+t^3)e_1,[e_2,e_4]\!=\!e_2+t^1e_3,[e_3,e_4]\!=\!t^2e_2+
(1+t^3)e_3.
\end{equation*}

Let us consider the picture including only the codifferentials $d_3$, $d_1^\sharp$
and $d_1(\lambda:\mu)$. The picture is very similar
to that of the moduli space of three dimensional Lie algebras.  The
family $d_1(\lambda:\mu)$ consists of a $\P^1$ with an action of the
symmetric group $\Sigma_2$, with orbifold points at $(1:1)$ and
$(1:-1)$.  The point $d_1(1:-1)$ has a\ jump deformation to
$d_3$, the four dimensional simple Lie algebra, while there is a
jump deformation from the point $d_1^\sharp$ to $d_1(1:1)$. The
point $d_1(1:0)$ is not an orbifold point, but is still special, with
a jump deformation to the point $d_2^\sharp$.  We did not see
this phenomenon in the 3 dimensional picture, but there was something
special about the point $d_2(1:0)$ in the family of codifferentials
$d_2(\lambda:\mu)$ (see Table \ref{Table 3}), because $\dim(H^1(d_2(1:0)))=2$,
instead of the generic value. Since $H^1$ influences the extensions of $\C$
by a Lie algebra, it is perhaps natural to expect that the 4 dimensional counterpart
to $d_2(1:0)$ should not behave generically.

\subsection{The codifferential $d_3(\lambda:\mu:\nu)$}

Before examining the cohomology in the generic case, we want to make
some general remarks about the family $d_3(\lambda:\mu:\nu)$, which we
will call \emph{the big family}, relating to the fact that the points
correspond to $\P^2/\Sigma_3$, in contrast to $d_3(\lambda:\mu)$,
 which we will refer to as
\emph{the small family}. Let us refer to elements in the orbit
of a point under the action of the symmetric group as conjugates. Most
points in $\P^2$ have precisely 6 conjugates, and the stabilizer of
the point is the trivial subgroup. The few exceptional cases are as
follows.
\begin{enumerate}
\item The points $(\lambda:\lambda:\mu)$, where $\lambda\ne\mu$ and
their conjugates are stabilized by a subgroup of order 2, so they each
have only 3 conjugates.
\item The point $(1:-1:0)$ and its 3 conjugates are also stabilized by
a subgroup of order 2.
\item The point $(1:r:r^2)$, where $r$ is a primitive cube root of
unity, and its 2 conjugates, are stabilized by the alternating group
$A_3$.
\item The point $(1:1:1)$ is stabilized by the entire group $\Sigma_3$.
\end{enumerate}
Next, consider the lines ($\P^1$) in $\P^2$ and the induced action of
$\Sigma_3$ on the set of lines. For most lines, the stabilizer of the
line is just the trivial subgroup, but again, there are a few
exceptions.
\begin{enumerate}
\item The line $(\lambda:\lambda:\mu)$ and its 3 conjugates are
stabilized by  subgroups of order 2.
\item The lines $(\lambda:\mu:c(\lambda+\mu))$ and their conjugates
are also stabilized by subgroups of order 2.
\end{enumerate}
It turns out that  when $c=0$ or $c=\pm1$, the cohomology of the
codifferentials corresponding to points on the line
$(\lambda:\mu:c(\lambda+\mu))$ does not follow the generic pattern.
The cohomology of the codifferentials corresponding to points on the
line $(\lambda:\lambda:\mu)$ is the same as the generic case with the
exception of the points $(1:1:0)$, $(1:1:2)$ and $(1:1:-2)$, which are
the points of intersection of this line with the three other special
lines.  Note also that the lines $(\lambda:\mu:c(\lambda+\mu))$ all
intersect in precisely the point $(1:-1:0)$, which makes this point
very special.

To determine the cohomology of a codifferential of type
$d_3(\lambda:\mu:\nu)$, read  Table \ref{Table 2} in descending order,
and whichever is the first pattern it matches, that gives its
cohomology.  However, we will present the description of the
cohomology in ascending order, because it is more natural to begin
with the generic pattern, and then proceed to the more exotic cases.

\subsection{The codifferential $d_3(\lambda:\mu:\nu)$: the generic case}
Generically, we have
\begin{align*}
H^1=&\langle\pha11(\lambda-\mu)+\pha21+\pha22(\mu-\nu)+\pha32,
\\&\pha11(-\lambda\mu+\lambda^2-\lambda\nu+\mu\nu)+\pha21(\lambda-\nu)+\pha31\rangle.
\end{align*}
For most generic values of $(\lambda:\mu:\nu)$ a natural basis to
choose for $H^2$ would be $H^2=\langle\psa{24}2,\psa{34}3\rangle$.
Then $d^\inf=d_3(\lambda,\mu+t^1,\nu+t^2)$, so there is no
difficulty in seeing what the deformations are equivalent to. However,
for certain generic values of the parameters, the two cocycles above are not
a basis of $H^2$, so we need to work with a more complex solution,
which yields a basis of $H^2$ for all generic values. Let us take
$$
H^2=\langle\psi_1=\psa{24}3,\psi_2=\psa{14}3\rangle.
$$
The universal infinitesimal deformation
$d^\inf=d_3(\lambda:\mu:\nu)+\psi_i t^i$ is miniversal, with base
$\A=\C[[t^1,t^2]]$. It is a bit more difficult to identify what the
miniversal deformation $d^\infty=d^\inf$ is equivalent to when we take
this more complicated basis of $H^2$.  In fact, if we let $x$ be a
root of the polynomial
$$
z^3+(-\nu+2\lambda-\mu)z^2+(\mu\nu-\lambda\nu-\lambda\mu+\lambda^2-t^1)z-t^2,
$$
and $y$ be a root of the polynomial
$$
z^2+(-x-\nu+\mu)z+x^2+x(\lambda-\mu)-t^1,
$$
then if $g$ is given by the matrix
$G=\left[\begin{smallmatrix}
1&0&0&0\\
x&1&0&0\\
x(\lambda-\mu)+x^2&y&1&0\\
0&0&0&1
\end{smallmatrix}\right]$, we have
$$
g^*(d^\infty)=d_3(\lambda+x:\mu+y-x:\nu-y).
$$
Thus for generic values of $(\lambda:\mu:\nu)$, all deformations of
$d_3(\lambda:\mu:\nu)$ simply move along this same big family.

\subsection{The codifferential $d_3(\lambda:\mu:-\lambda-\mu)$}
In this case $H^1$ is the same as the generic case,
and for most values of $\mu$, one can use the cocycles $\psa{14}1$ and
$\psa{24}2$ as the basis of $H^2$. Since the brackets of these two
cocycles vanish, the resulting infinitesimal deformation
$$d^\inf=d_3(\lambda:\mu:-\lambda-\mu)+\psa{14}1t^1+\psa{24}2t^2$$ is
miniversal, and in fact coincides with
$d_3(\lambda+t^1,\mu+t^2,-\mu-\lambda)$. In this case, it is obvious
that the deformations of our codifferential just lie along the big
family. The values for which these two elements do not form a basis of
$H^2$ are $(1:-1:0)$, which we will cover separately, and $(1:1:-2)$,
for which $\psi_1=\psa{14}3$ and $\psi_2=\psa{24}2$ give a
basis of $H^2$. It is still true that the brackets of these cocycles
vanish, and the deformations lie along the big family, although the
expression of the member of the family corresponding to the element
$d^\infty$ is more complicated in this case, and will be omitted.

Thus the family $d_3(\lambda:\mu:-\lambda-\mu)$ is not special in
deformation theory.  This is a bit surprising, since $H^3$ does not
vanish for elements of this subfamily, so it would not have been unreasonable
to expect
that there would be some obstructions to the extension of an
infinitesimal deformation.

\subsection{The codifferential
$d_3(\lambda:\mu:0)=\mathfrak{r}_{3,\mu/\lambda}(\C)\oplus\C$}.
The dimensions of $H^1$ and $H^2$ increase to 3,
and $H^3$ is 1-dimensional as well. The two cocycles $\psi_1$ and
$\psi_2$ chosen as basis elements for $H^2$ in the generic case remain nontrivial and one can
find an independent nontrivial cocycle $\psa{13}1+\psa{23}2$.
However, this choice of a basis turns out to be inconvenient, and a
slight modification of the basis will make the presentation simpler.
We have
\begin{align*}
H^1=&\langle\pha11\lambda(\lambda-\mu)+\pha21\lambda+\pha31,
\pha11+\pha22+\pha33,\pha43\rangle\\
H^2=&\langle\psi_1=\psa{24}3+\psa{14}3\lambda,\psi_2=\psa{14}3,\psi_3=\psa{13}1+\psa{23}2\rangle.
\end{align*}

Let $d^\inf=d_3(\lambda:\mu:0)+\psi_it^i$.
We compute
\begin{align*}
[\psi_1,\psi_3]&=-\pha{124}1+\pha{234}3+\pha{134}3\lambda+\pha{124}2\lambda=-D(\psa{12}2+\psa{13}3)\\
[\psi_2,\psi_3]&=\pha{134}3+\pha{124}2,
\end{align*}
so that
$$
\tfrac12[d^\inf,d^\inf]=-D(\zeta_1)t^1t^3+(\pha{134}3+\pha{124}2) t^2t^3,
$$
where $\zeta_1=\psa{12}2+\psa{13}3$.

Note that in the case $t^3=0$, since $\psi_1$ and $\psi_2$ span
the same subspace as the ones we used in the generic case, a
deformation with $t^3=0$ is equivalent to one in the family. Thus
there is a two parameter family of deformations of
$d_3(\lambda:\mu:0)$ along the big family $d_3(\alpha:\beta:\eta)$.

On the other hand,  when $t^3$ does not vanish, we will have to
consider how $d^\inf$ extends to a higher order deformation. It turns out
that when $\lambda=\mu$, the codifferential $\pha{134}3+\pha{124}2$ is
a coboundary, but otherwise, it can be taken as a basis of $H^3$, and
so is an obstruction to the extension of $d^\inf$ to a higher order
deformation.  We will first consider this obstructed case.

\subsubsection{$\mu\ne\lambda$} In this case we have
$$
H^3=\langle\phi=\pha{134}3+\pha{124}2\rangle.
$$
We extend $d^\inf$ to the second order deformation
$$
d^2=d^\inf+\zeta_1t^1t^3.
$$
Since $\tfrac12[d^2,d^2]=\phi t^2t^3$, the second order deformation is
miniversal and the base of the miniversal deformation is
$\A=\C[[t^1,t^2,t^3]]/(t^2t^3)$.

Thus any true deformation is given by taking $d^\infty=d^2$ with
either $t^2=0$ or $t^3=0$.  Since the case $t^3=0$ has already been
examined, we consider the case $t^2=0$. In this case,  the matrix $A$
of the deformation $d^\infty$ is given by
$A=\left[\begin{smallmatrix}
0&t^2&0&\lambda&1&0\\
t^1t^3&0&t^3&0&\mu&1\\
0&t^1t^3&0&\lambda t^1&t^1&0\\
0&0&0&0&0&0
\end{smallmatrix}\right]$. When $t^1=-\tfrac{(\lambda-\mu)^2}4$,
then
$d^\infty\sim d_1(1:0)$.  Since $\lambda\ne\mu$, note that this
deformation is not a jump deformation of $d_3(\lambda:\mu:0)$ but
occurs some ``distance" away from the original codifferential.

When $t^1\ne-\tfrac{(\lambda-\mu)^2}4$,  then 
$d^\infty=d_2^\sharp$. This is a jump deformation, since it is
independent of the value of $t^1$, as long as it is small.

Thus we obtain that the deformations of $d_3(\lambda:\mu:0)$,  for
$\lambda\ne\mu$, live along two planes in the $(t^1,t^2,t^3)$ space.
One is the plane $t^3=0$ determining deformations along the family
$d_3(\lambda:\mu:\nu)$, while those which lie in the plane $t^2=0$
are  equivalent to $d_2^\sharp$, except along the line
$t^1=-\frac{(\lambda-\mu)^2}4$, which is not important to us, because
this line does not include the origin. We say that a family of
deformations is \emph{not local} if the origin in the $t$-parameter
space is not part of the family. Thus the deformations along the line
$t^1=-\frac{(\lambda-\mu)^2}4$ are not local, in this sense. Only
local deformations play a role in determining how the moduli space is
glued together.

\subsubsection{$\mu=\lambda$}
This is the codifferential $d_3(1:1:0)=\mathfrak{r}_3(\C)\oplus\C$. We have
$$
[\psi_2,\psi_3]=\pha{134}3+\pha{124}2=-D(\zeta_2),
$$
where $\zeta_2=-\psa{23}3+\psa{13}3+\psa{24}4+\psa{12}2$. So the
second order deformation is given by
$$
d^2=d^\inf+\zeta_1t^1t^3+\zeta_2t^2t^3.
$$
Since $\ph=\pha{124}3$ is a nontrivial 3-cocycle, we can take
$$
H^3=\langle\phi=\pha{124}3\rangle.
$$
Now
$$
\tfrac12[d^2,d^2]=-2\phi t^2t^3(t^1+t^3)-D(\zeta_3)t^2(t^3)^2+
(\pha{124}4-\pha{123}3)t^2(t^3)^2(t^1+t^2),
$$
where $\zeta_3=\psa{23}4$.  We did not obtain any second order relations, but
because of the term involving $\ph$ in the bracket above,
there is a third order relation $t^2t^3(t^1+t^2)$.
The last term in the bracket is of higher order, so can be ignored in computing the
third order deformation. We can take
$$
d^3=d(1:1:0)+\psi_it^i +\zeta_1t^1t^3+\zeta_2t^2t^3+\zeta_3t^2(t^3)^2.
$$
One computes that
\begin{align*}
\tfrac12[d^3,d^3]=
&-2\phi t^2t^3(t^1+t^3)+2(\pha{124}4-\pha{123}3)t^2(t^3)^2(t^1+t^2)\\
&+2\pha{123}4t^2(t^3)^3(t^1+t^2).
\end{align*}
But this term is equal to zero, using the third order relation.  Thus
the base of a versal deformation is
$\A=\C[[t^1,t^2,t^3]]/(t^2t^3(t^1+t^2))$, $d^\infty=d^3$, and the
formal deformation corresponds to an actual deformation along the
three planes $t^3=0$, $t^2=0$ and $t^2=-t^1$.

The plane $t^3=0$ corresponds to the generic case, which gives a 2
parameter space of deformations along the family
$d_3(\lambda:\mu:\nu)$.

Now consider the plane spanned by $t^2=0$. If neither $t^1$ nor $t^3$
vanish, we then have $d^\infty\sim d_2^\sharp$, so there is a jump
deformation from $d_3(1:1:0)$ to $d_2^\sharp$, just as for the other
points on the line $d_3(\lambda:\mu:0)$.

When $t^2=t^3=0$, we are on the plane $t^3=0$, which we discussed
already.  When $t^2=t^1=0$, then
we get a jump deformation to the point $d_1(1:0)$.  This is not like
the generic case.

Finally, let us consider the case when $t^2=-t^1$ and $t^3\ne0$. Let
us express $t^1=\tfrac{\alpha\beta}{(\alpha+\beta)^2}$, and let
$x=\tfrac{\alpha+\beta}{t^3}$. Then $g^*(d^3)=d_1(\alpha,\beta)$ if
$g$ is given by the matrix
$G=\left[\begin{smallmatrix}
-x&0&0&0\\
0&-x&0&0\\
0&\frac{\beta}{t^3}&-\frac1{t^3}&x\\
0&-\beta&1&0
\end{smallmatrix}\right].$
Note that $(\alpha:\beta)$ is independent of $t^3$ as long as $t^3$
is nonzero. On the plane $t^2=-t^1$, we see that for $t^1=0$, the
deformation jumps to $d_1(1:0)$, but when $t^1\ne0$, we deform along
the family $d_1(\alpha:\beta)$. It is as if $d_3(1:1:0)$ sits just
above the point $d_1(1:0)$ and deforms as if it were that element.
This is the usual pattern we have already discussed when there is a
jump deformation to a point.

\subsection{The codifferential $d_3(\lambda:\mu:\lambda+\mu)$}
We exclude the codifferential $d_3(1:0:1)$ from our consideration
here, because it coincides with the codifferential $d_3(1:1:0)$, which
we treated previously. $H^1$ is the same as the generic case, and we have
\begin{align*}
H^2=\langle\psi_1=&2\psa{12}1\lambda\mu(\lambda+\mu)^2+
\psa{12}2\lambda\mu^2(\lambda+\mu)(\lambda+2\mu)\\
  &+\psa{12}3\lambda^2\mu^2(\lambda+\mu)^2-
  \psa{13}1\mu(\mu^2+2\lambda\mu+2\lambda^2)\\
&-\psa{13}2\mu^2(\lambda+\mu)+\psa{23}1(\mu^2+\lambda\mu+\mu^2)\\
  &+(\psa{14}1+\psa{23}2)\lambda^2\mu^2(\lambda+\mu)\\
\psi_2=&\psa{14}1\\
\psi_3=&\psa{14}2\rangle.
\end{align*}
The brackets of $\psi_1$ with itself and $\psi_3$ are coboundaries,
its bracket with $\psi_2$ is a nontrivial cocycle, and the rest of the
brackets vanish. From this, one sees immediately that the second order
relation is $t^1t^2=0$, but it is not so obvious what higher order
terms might be necessary to add in order to obtain the relation on the
base of the miniversal deformation. Since the space of 3-cocycles is
12 dimensional, we know that a miniversal deformation can be
expressed in the form
\begin{equation*}
d^\infty=d_3(\lambda,\mu,\lambda+\mu)+\psi_i t^i +\zeta_i x^i,
\end{equation*}
where $\zeta_1,\dots,\zeta_{11}$ is a pre-basis of the space of
3-coboundaries. In fact, we can give this pre-basis as
$$
\{\zeta_i,i=1,\dots11\}=
\{\psa{12}1,\psa{12}2,\psa{12}3,\psa{12}4,\psa{13}1,\psa{13}2,
\psa{13}3,\psa{13}4,\psa{23}1,\psa{23}2,\psa{23}4
\}.
$$
Note that the first 10 of these vectors are just the first 10
elementary 2-cochains. Also
$$
H^3=\langle\phi=\pha{124}3\rangle,
$$
and we can complete the linearly independent set given by the
$D(\zeta_i)$ and $\phi$ to a basis
$\{D(\zeta_1),\dots,D(\zeta_11),\phi,\tau_1,\dots\tau_4\}$ of $L_3$.
Then we must have
$$
[d^\infty,d^\infty]=D(\zeta_1)s^1+\cdots+\D(\zeta_{11})s^{11}+\phi s^{12}
  +\tau_1 s^{13}+\cdots \tau_4 s^{16},
$$
for some coefficients $s^1,\dots s^{16}$, where these coefficients are
expressed as polynomials in the variables $t^i$ and $x^i$. Now all of
these coefficients must be equal to zero, once you take into account
the relation on the base of the miniversal deformation, which is the
coefficient $s^{12}$.  The expression one obtains for $s^{12}$ by
direct computation from the form of $d^\infty$ will have the variables
$x^i$ in it, but it should depend only on the variables $t^i$.  The
trick is to solve the first 11 equations for $x^i$ as functions of the
variables $t^i$, and then substitute these into the formula for
$s^{12}$ to obtain the relation on the base.

The relation on the base of the miniversal deformation is simply
$t^1t^2=0$,  which is exactly the second order relation. If you solve
for the coefficients of $s^{13},\dots,s^{16}$, then they turn out to
be multiples of $s^{12}$, so they are equal to zero using the relation
on the base.

Let us study the deformations of
$d(\lambda:\mu:\lambda+\mu)$. Since the relation on the base of the
miniversal deformation is $t^1t^2=0$, in any true deformation, we must
have either $t^1=0$ or $t^2=0$.

When $t^1=0$, then $d^\infty=d_3(\lambda:\mu:\lambda+\mu)+\psi_2
t^2+\psi_3 t^3$, so that for any values of $t^2$ and $t^3$ we have a
deformation along the big family.  In fact, $d^\infty\sim
d_3(\alpha:\beta:\eta)$ where
\begin{align*}
\alpha=&\tfrac12(\lambda+\mu+t^2+\sqrt{(t^2+\lambda-\mu)^2+4t^3})\\
\beta=&\tfrac12(\lambda+\mu+t^2-\sqrt{(t^2+\lambda-\mu)^2+4t^3})\\
\eta=&\lambda+\mu.
\end{align*}

The interesting case is when $t^2=0$. The matrix of $d^\infty$ is
quite complicated, so we won't reproduce it here, but it should be
noted that some terms have $t^3-\lambda\mu$ in the denominator, so
that $t^3=\lambda\mu$ may not correspond to an actual deformation.
When $t^1\ne0$, then
$d^\infty\sim\mathfrak g_8\left(\tfrac{\lambda\mu-t^3}{(\lambda+\mu)^2}\right)$.
In particular, if we set $t^3=0$, we see that
there is a jump deformation to $d_1(\lambda:\mu)$, and that we also deform along the family
$d_1(\lambda:\mu)$ when $t^3\ne 0$.

\subsection{The codifferential $d_3(1:-1:0)$}
For this codifferential, from the fact that $H^2$ and $H^3$ both have
dimension 5, we expect to see some interesting phenomena, both because
the tangent space to the space of deformations has dimension 5, and
since $H^3$ has high dimension, the dimension of the variety of
deformations would likely be lower than 5. We can give bases for the cohomology as follows:
\begin{align*}
H^1=\langle&2\pha11+\pha22+\pha31,\pha11+\pha22+\pha33,\pha43\rangle\\
H^2=\langle&\psi_1=\psa{24},\psi_2=\psa{14}3,\psi_3=\psa{12}3-
  \psa{13}3-\psa{23}3+\psa{14}4+\psa{24}4\\
&\psi_4=\psa{23}1-2\psa{23}2,\psi_5=\psa{12}4-\psa{13}4-\psa{23}4\rangle\\
H^3=\langle&\phi_1=\pha{124}1,\phi_2=\pha{124}3,\phi_3=\pha{123}2+
  \pha{123}3-\pha{234}4\\
&\phi_4=\pha{123}4,\phi_5=\pha{123}2+\pha{124}4-\pha{234}4\rangle.
\end{align*}
A pre-basis of the 3-coboundaries is
$$
\{\zeta_1,\dots,\zeta_9\}=
\{\psa{12}1,\psa{12}2,\psa{12}3,\psa{12}4,\psa{13}1,\psa{13}2,\psa{13}3,
\psa{13}4,\psa{14}4\}.
$$

A miniversal deformation is given by
$$
\d^\infty=d_3(1:-1:0)+\psi_i t^i+\zeta_i x^i,
$$
where the $x^i$ are expressible as power series in the variables
$t^i$. Since not all of the brackets of the $\psi_i$ vanish, we do not
expect that the coefficients $x^i$ are all equal to zero, in general.

We can express
$$
[d^\infty,\d^\infty]=D(\zeta_i)s^i +\phi_i s^{9+i} +\tau_i s^{14+i},
$$
where $D(\zeta_i)$, $\phi_i$ and $\tau_i$ form a basis of $L_3$.
Solving $s^1=\cdots=s^9=0$ for $x^1,\dots,x^9$  in terms of
$t^1,\dots, t^9$, and substituting these values of the $x^i$ into the
formulas for $s^{10},\dots,s^{14}$, we obtain 5 relations on the base of the
versal deformation, the
simplest of which is
$$
\tfrac{(t^3(t^1)^2+4t^1t^2t^4+4t^2t^4)}{t^1-2}=0.
$$
Some of these relations have $t^1-1$ or $t^1-2$ as a factor of the
denominator, which means that there may not be a solution when $t^1$
takes on these values. There should be an actual, rather than just a
formal power series solution for all  values of $t^i$ which make all 5
of the relations vanish. When we solved for the zeros of the
relations,  we obtained the following 5 solutions:
\begin{align*}
1)\quad t^1=&t^2=t^4=0\\\
2)\quad t^3=&t^4=t^5=0\\
3)\quad t^2=&t^3=t^5=0\\
4)\quad t^3=&t^5=0,\qquad t^1=-1\\
5)\quad t^2=&\tfrac{-t^1(t^1-2)^2}8,t^3=\tfrac{t^4(t^1-2)^2(t^1+1)}{2t^1},
t^5=\tfrac{(t^4)^2(t^1-2)^2(t^1+1)}{(t^1)^2}.\\
\end{align*}
Note that each of these solutions is only a 2-dimensional subvariety
of the 5-dimensional tangent space.

For the first solution,  the matrix of the corresponding $d^\infty$ is
$$
A=\left[\begin{smallmatrix}
0&0&0&1&1&0\\
0&0&0&0&-1&1\\
t^3&-t^3&-t^3&0&0&0\\
t^5&-t^5&-t^5&t^3&t^3&0
\end{smallmatrix}\right].
$$
Along the curve $t^5=(t^3)^2$, $d^\infty\sim d_1(1:-1)$.  For all
other points on the $(t^3,t^5)$-plane, $d^\infty\sim d_3$.  This fits
with our prior observation that there is a jump deformation from
$d_1(1:-1)$ to $d_3$.  Thus we have jump deformations from
$d_3(1:-1:0)$ to both $d_1(1:-1)$ and $d_3$.

For the second solution,  the matrix of $d^\infty$ is
$$
A=\left[\begin{smallmatrix}
0&0&0&1&1&0\\
0&0&0&0&-1+t^1&1\\
0&0&0&t^2&0&0\\
0&0&0&0&0&0
\end{smallmatrix}\right].
$$
Deformations corresponding to this solution give a two parameter
family of deformations along the big family $d_3(\lambda:\mu:\nu)$.

For the third solution,  the matrix of $d^\infty$ is
$$
A=\left[\begin{smallmatrix}
0&-t^4t^1&t^4&1&1&0\\
0&0&-2t^4&0&-1+t^1&1\\
0&0&0&0&0&0\\
0&0&0&0&0&0
\end{smallmatrix}\right].
$$
When $t^4=0$, this is just a special case of the previous solution,
and in fact, in this case $d^\infty=d_3(1:-1+t^1:0)$. Supposing that
$t^4\ne0$, then when $t^1\ne-1$, then $d^\infty\sim\d_2^\sharp$ which
is a jump deformation.

The fourth solution has $t^1=-1$, which means it is not local, so does
not contribute to our picture of the moduli space. Although the
solution is interesting, we will omit it here.

For the last solution, which is the most complicated of them all, the first three columns of the matrix of $d^\infty$ are
$$
\left[\begin{smallmatrix}
0&t^1t^4&t^4\\
\frac{t^4(t^1-2)^2(t^1+1)}4&0&-2t^4\\
\frac{t^4(t^1-2)^4(t^1+1)}{8t^1}&\frac{-t^4(t^1-2)^2(-4+3(t^1)^2)}{8t^1}&
\frac{-t^4(t^1-2)^2(t^1+1)}{2t^1}\\
\frac{(t^4)^2(t^1+2)(t^1-2)^4(t^1+1)}{8(t^1)^2}&\frac{-(t^4)^2((t^1)^3-2(t^1)^2+4t^1+8)(t^1-2)^2}{8(t^1)^2}&
\frac{-(t^4)^2(t^1+1)(t^1-2)^2}{(t^1)^2}\\
\end{smallmatrix}\right].
$$
Note that $t^1$ appears in the denominator, so cannot vanish for this
solution.  If $t^4=0$, $t^1=-1$ or $t^1=2$, then the fifth solution
coincides with one of the previous four, so we will not consider these
cases here.
The matrix of $A$ is so complicated that in order to determine
which standard form the codifferential is equivalent to
three, we first had to transform $A$ into a matrix of an equivalent
codifferential which had a simpler matrix.  We found that
$d^\infty\sim d_1(\alpha:\beta)$ where
$$
\alpha=\tfrac{t^1+\sqrt{5(t^1)^2-16t^1+16}}2,\quad
\beta=\tfrac{t^1-\sqrt{5(t^1)^2-16t^1+16}}2.
$$
Note that if we were to set $t^1=0$ in the above, we would obtain the
codifferential $d_1(1:-1)$, to which we already obtained a jump
deformation in the first solution above.

The picture of the local deformations of $d_3(1:-1:0)$ is as follows.
First, we can deform along the big family. Secondly, we can deform to
$d_2^\sharp$, like any other member of the family
$d_3(\lambda:\mu:0)$. Thirdly, like any other member of the family
$d_3(\lambda:\mu:\lambda+\mu)$, we have a jump deformation to an
element in the family $d_1(\lambda:\mu)$. Because the element we
deform to is $d_1(1:-1)$, which has an extra deformation to the
element $d_3$, we can  also deform to this element, as well as deforming
along the family $d_1(\lambda:\mu)$.

\subsection{The codifferential $d_3(\lambda:\mu)$}
This family does not have an action of the symmetric group, which is
important to keep in mind.
Generically, $H^1$ and $H^2$ are 4 dimensional.  The generic basis of $H^2$ below
consists of elements which are linearly independent nontrivial
cocycles for generic values of $\lambda$ and $\mu$ except in the
special case $\lambda=\mu$, which we will treat separately.  Of
course, for those values of $\lambda$ and $\mu$ for which $\dim
H^2>4$, they do not span $H^2$. Generically, we have
\begin{align*}
H^1=\langle&\psa12,\psa22(\lambda-\mu)+\psa32,\psa22+\psa33,\psa21(\lambda-\mu)+\psa31\rangle\\
H^2=\langle&\psi_1=\psa{34}3,
\psi_2=\psa{14}2,
\psi_3=\psa{24}1,
\psi_4=\psa{24}2
\rangle.
\end{align*}
All of the brackets of these nontrivial cocycles vanish, so the
miniversal deformation is just the first order deformation
$d^\infty=d_3(\lambda:\mu)+\psi_i t^i$, and there are no relations on
the base. The matrix of $d^\infty$ is given by
$A=\left[\begin{smallmatrix}
0&0&0&\lambda&t^3&0\\
0&0&0&t^2&\lambda+t^4&1\\
0&0&0&0&0&\mu+t^1\\
0&0&0&0&0&0
\end{smallmatrix}\right]$.

If $t^3\ne 0$, then $d^\infty\sim d_3(\alpha,\beta,\eta)$, where
\begin{equation*}
\alpha=\lambda+\tfrac{t^4+\sqrt{(t^4)^2+4t^2t^3}}2,
\beta=\lambda+\tfrac{t^4-\sqrt{(t^4)^2+4t^2t^3}}2,\eta=\mu+t^1.
\end{equation*}
If $t^1=t^2=t^4=0$, then  $d^\infty\sim d_3(\lambda:\lambda:\mu)$,
so there is a jump deformation from $d_3(\lambda:\mu)$ to
$d_3(\lambda:\lambda:\mu)$. Thus we see that the codifferential
$d_3(\lambda:\mu)$ sits over the codifferential
$d_3(\lambda:\lambda:\mu)$ and deforms along the big family as if it
were that codifferential. All of the deformations of
$d_3(\lambda:\mu)$ which do not lie along the hyperplane $t^3=0$ lie
along the big family.

Now, consider the hyperplane $t^3=0$.  The eigenvalues of the submatrix
$B=\left[\begin{smallmatrix}
\lambda&0&0\\
t^2&\lambda+t^4&1\\
0&\mu+t^1
\end{smallmatrix}\right]$ are $\lambda$, $\lambda+t^4$ and $\mu+t^1$.
If these eigenvalues are all distinct, then $d^\infty\sim
d_3(\lambda:\lambda+t^4:\mu+t^1)$.  Otherwise, one of the conditions
$t^4=0$, $t^1=\lambda-\mu$, or $t^1-t^4=\lambda-\mu$ holds.  Of these
conditions, only the first one is local, so we will not consider the
other two.
%
Consider the plane $t^3=t^4=0$. Unless $t^2=0$ or $t^1=\lambda-\mu$,
$d^\infty$ is still equivalent to $d_3(\lambda:\lambda+t^4:\mu+t^1)$.
Again, the second condition is not local, so we will ignore it. On the
line $t^2=0$, we have  $d^\infty\sim d_3(\lambda:\mu+t^1)$, so we get
a deformation along the $d_3(\lambda:\mu)$ family.

%
To summarize the generic deformation behavior of an element of the
family $d_3(\lambda:\mu)$, we have the following picture. First, we
can always deform along the family to which an element belongs, so
there is a deformation along the family $d_3(\lambda:\mu)$. Secondly,
there is a jump deformation to the element $d_3(\lambda:\lambda:\mu)$
in the big family. Whenever there is a jump deformation, then we can
deform in any manner in which the element we jump to deforms, and thus
there is a deformation along the big family as well.  Note that the
line $(\lambda:\lambda:\mu)$, which is one of the lines in $\P^2$
with nontrivial stabilizer, is the target of our jump deformations, so
the elements $d_3(\lambda:\lambda:\mu)$ are special not in the sense
that they have more deformations, but that there are extra
deformations to them.

\subsection{The codifferential $d_3(1:1)$}
Even though the the dimension of $H^2$ for this element is the same as
the generic case of $d_3(\lambda:\mu)$, we have to use a different
basis for $H^2$ than in the generic case.
\begin{equation*}
H^2=\langle\psi_1=\psa{14}3,\psi_2=\psa{14}1,\psi_3=\psa{24}1,\psi_4=
\psa{24}3\rangle.
\end{equation*}
As in the generic case, the brackets of these cocycles all vanish, so
the universal infinitesimal deformation is the miniversal deformation
$d^\infty$, with matrix
$
A=\left[\begin{smallmatrix}
0&0&0&1+t^2&t^3&0\\
0&0&0&0&1&1\\
0&0&0&t^1&t^4&1\\
0&0&0&0&0&0
\end{smallmatrix}\right]
$.
When $t^1\ne0$, we obtain a complicated deformation along the family
$d_3(\alpha:\beta:\eta)$. To understand the solution a bit better,
when we solve for a matrix transforming $A$ into one representing a
codifferential of the
form $d_3(\alpha:\beta:\eta)$, we obtain a solution which satisfies
\begin{align*}
t^3=&\tfrac{1/3(\alpha^3+\beta^3+\eta^3)(t^2+3)^3+(\alpha+\beta+
\eta)(t^2+3)p(\alpha,\beta,\eta,t^2)+8(\alpha+\beta+\eta)^3}
{t^1(\alpha+\beta+\eta)^3}\\
t^4=&\tfrac{1/3(\alpha^2+\beta^2+\eta^2)(t^2+3)^2-1/2(\alpha+\beta+
\eta)((t^2)^2+2t^2+3))}{(\alpha+\beta+\eta)^2},
\end{align*}
where $p$ is a polynomial which is homogeneous, quadratic and
symmetric in $\alpha$, $\beta$ and $\eta$ and quadratic in $t^2$.
Consequently,  when $t^2\ne-3$, for any values of $t^2$, $t^4$ and
$t^3$, we obtain exactly one solution up to the action of the
symmetric group, and thus one member of the family
$d_3(\alpha:\beta:\eta)$ is determined. This follows since the line
$\alpha+\beta+\eta=0$ intersects the quadric surface determined by the
equation for $t^4$ above in exactly the orbifold points
$(1:\tfrac{-1+\sqrt3}2:\tfrac{-1-\sqrt3}2)$ and
$(1:\tfrac{-1-\sqrt3}2:\tfrac{-1:\sqrt3}2)$, which do not lie on the
cubic surface determined by the equation for $t^3$.

When $t^1\ne0$ and $t^2=t^3=t^4=0$, then $d^\infty\sim d_3(1:1:1)$, so
there is a jump deformation to this element, as we expect from the
generic case.


When $t^1=0$, then the eigenvalues of the submatrix
$B=\left[\begin{smallmatrix}
1+t^2&t^3&0\\
0&1&1\\
t^1&t^4&1
\end{smallmatrix}\right]$ are $1+t^2$ and $1\pm\sqrt{t^4}$, so they
are distinct unless $t^4=0$ or $t^4=(t^2)^2$. Thus, except in these
two cases we have $d^\infty\sim d_3(1+t^2,1+\sqrt{t^4},1-\sqrt{t^4})$.
On the plane $t^1=t^4=0$ we have $d^\infty\sim d_3(1+t^2,1,1)$.

On the surface $t^1=0$, $t^4=(t^2)^2$ except on the curve $t^3=0$ we
have $d^\infty\sim d_3(1+t^2,1+t^2,1-t^2)$. Finally, on the curve
$t^3=0$ on this surface we have $d^\infty\sim d_3(1+t^2,1-t^2)$, so we
obtain a deformation along the family $d_3(\lambda:\mu)$ on this
curve.

Thus, just like any other generic value, there is one curve along
which there is a jump deformation to the corresponding point
$d_3(1:1:1)$ on the large family, another curve along which we deform
along the $d_3(\lambda:\mu)$ family, and otherwise,  all deformations
are along the big family.  In a way, it is surprising that the one
point in $\P^2$ which is fixed by every permutation does not have any
special properties in terms of deformation theory, but as we have
seen, there just isn't anything particularly special about the
deformations of this codifferential.

\subsection{The codifferential $d_3(1:-2)$}
We have
\begin{align*}
H^2=&\langle\psi_1=\psa{14}2,\psi_2=\psa{24}3,\psi_3=\psa{34}3,\psi_4=
\psa{24}1\rangle\\
H^3=&\langle\pha{123}4\rangle.
\end{align*}
Even though $H^3\ne0$, it turns out that the brackets of all the
$\psi$'s with each other vanish, so the miniversal deformation
$d^\infty$ coincides with the infinitesimal deformation, and its
matrix is given by $
A=\left[\begin{smallmatrix}
0&0&0&1&t^4&0\\
0&0&0&t^1&1&1\\
0&0&0&0&t^2&-2+t^3\\
0&0&0&0&0&0
\end{smallmatrix}\right]$.
Because this matrix has no terms on the left hand
side, it is natural to guess that the deformations are either along the family
$d_3(\alpha:\beta:\eta)$ or the family $d_3(\lambda:\mu)$, with
possibly a few exceptional codifferentials.

When $t^1\ne 0$ we have a solution of the form
\begin{align*}
t^3=&\tfrac{\alpha+\beta+\eta}{q}\\
t^2=&-\tfrac{(\alpha+\eta-2q)(\alpha+\beta-2q)(\beta+\eta-2q)}{q^2(\alpha+
\beta+\eta-3q)}\\
t^4=&\tfrac{-(\alpha-q)(\beta-q)(\eta-q)}{t^1q^2(\alpha+\beta+\eta-3q)},
\end{align*}
where $q$ is a nonzero free parameter. These equations are symmetric
in $\alpha$, $\beta$ and $\eta$. If $t^3\ne0$, then
$\alpha+\beta+\eta\ne 0$, and we can solve the first equation for
$q$ and get

\begin{align*}
t^2=&\tfrac{-((\alpha+\beta)t^3-2(\alpha+\beta+\eta))
((\beta+\eta)t^3-2(\alpha+\beta+\eta))((\alpha+\eta)t^3-2(\alpha+\beta+\eta)}{(\alpha+\beta+\eta)^3(t^3-3)}\\
t^4=&\tfrac{-(\alpha t^3-(\alpha+\beta+\eta))(\beta t^3-(\alpha+\beta+\eta))(\eta t^3-(\alpha+\beta+\eta))}
{t^1(t^3-3)(\alpha+\beta+\eta)^3}.
\end{align*}
 We can express these equations in the form
\begin{align*}
t^2=&\tfrac{\alpha\beta\eta (t^3)^3 +(\alpha+\beta+\eta)(t^3-2)p(\alpha,\beta,\eta,t^3)}
{(\alpha+\beta+\eta)^3(t^3-3)}\\
t^4=&\tfrac{-\alpha\beta\eta (t^3)^3 +(\alpha+\beta+\eta)r(\alpha,\beta,\eta,t^3)}
{t^1(t^3-3)(\alpha+\beta+\eta)^3},
\end{align*}
where $p$ and $q$ are homogeneous, quadratic and symmetric in
$\alpha$, $\beta$ and $\eta$. The surfaces represented by these two
equations are both cubic, so there are 9 points of intersection.
Since every cubic which is given by a symmetric, homogeneous
polynomial either contains the line $\alpha+\beta+\eta=0$ or
intersects this line in precisely the points $(1:-1:0)$, $(1:0:-1)$
and $(0:1:-1)$, there are six points in the intersection of these two
cubics not on this line, which uniquely determine the  codifferential
$d_3(\alpha:\beta:\eta)$ to which $d^\infty$ is equivalent. The matrix
representing the transformation can be chosen with nonzero
determinant, as long as $t^3\ne3$.  The condition $t^3\ne0$ can also
be overcome, because if we substitute $t^3=0$ in the above, then the
problem still has a solution.Thus, whenever $t^1\ne0$  and $t^3\ne 3$,
the deformation is equivalent to a member of the family
$d_3(\alpha:\beta:\eta)$.

When $t^1=0$, then as long as $t^2\ne0$ and $t^3\ne1$, $d^\infty\sim
d_3(\alpha:\beta:1)$, where
$\alpha=\frac{t^3-1+\sqrt{(t^3-3)^2+4t^2}}2$ and
$\beta=\frac{t^3-1-\sqrt{(t^3-3)^2+4t^2}}2$.

When $t^1=t^2=0$ and $t^4\ne 0$ we have $d^\infty\sim d_3(1:1:t^3-2)$.
As a consequence, if we set $t^3=0$, we have a jump deformation from
$d_3(1:-2)$ to $d_3(1:1:-2)$. On the other hand, when $t^1=t^2=t^4=0$,
then  $\d^\infty\sim d_3(1:t^3-2)$. When $t^1=0$ and $t^3=1$, then we
also have a deformation along the big family.
%
The upshot of all this analysis is that $d_3(1:-2)$ is really not
special in terms of deformation theory. It deforms along its own
family, jumps to $d_3(1:1:-2)$, and deforms along that family.

\subsection{The codifferential $d_3(1:2)$}
We have
\begin{align*}
H^2=\langle&\psi_1=\psa{12}2-\psa{12}3+\psa{13}2+\psa{13}3,
\psi_2=\psa{14}1+\psa{34}3,\\&\psi_3=\psa{14}2,\psi_4=\psa{34}3,\psi_5=
\psa{24}1\rangle\\
H^3=\langle&\pha{124}3+\pha{134}2\rangle.
\end{align*}
This time, we do have some nonzero brackets, but only those brackets
of $\psi_1$ with $\psi_2$, $\psi_4$ and $\psi_5$, with the first one
being a nontrivial cocycle, so that the second order relation is
$t^1t^2=0$.  After some work, one obtains that there is one relation on the
base of the versal deformation,
$$
t^1t^2(-1-t^4+t^3t^5)=0,
$$
so that there are three distinct solutions for a true deformation,
given by the three factors of the miniversal deformation. Notice that
the third factor does not give rise to a local deformation.

Let us study the first solution,  when $t^1=0$. This case is simplest.
The matrix corresponding to $d^\infty$ is
$
A=\left[\begin{smallmatrix}
0&0&0&1+t^4&t^5&0\\
0&0&0&t^3&1&1\\
0&0&0&0&0&2+t^2+t^4\\
0&0&0&0&0&0
\end{smallmatrix}\right]$.
%

When $t^5\ne0$, then $d^\infty$ is equivalent to
$d_3(\alpha:\beta:2+t^2+t^4)$, where
$\{\alpha,\beta\}=\tfrac{2+t^4\pm\sqrt{(t^4)^2+4t^3t^5}}2$. If
$t^2=t^3=t^4=0$, then the deformation is equivalent to $d_3(1:1:2)$
for all $t^5\ne0$, giving the expected jump deformation.

What happens if $t^5=0$? As long as $t^2\ne-1$ and
$t^2t^3+t^4t^3+t^3+t^4\ne0$, then the deformation is still along the
big family. 
 If $t^2\ne -1$, but
$t^2t^3+t^4t^3+t^3+t^4=0$, then as long as $t^3\ne0$ and $t^4\ne0$, the
deformation is in the big family.  If $t^4=0$, then $t^3=0$ or
$t^2=-1$, and in both cases we deform along the family
$d_3(\alpha:\beta)$. Thus, the first solution to the relations on the
base does not have any surprises.

The second solution to the relations on the base is $t^2=0$. We may as
well assume that $t^1\ne 0$ and that $-1-t^4+t^3t^5\ne0$ for this
case. Then the matrix of $d^\infty$ is
$
A=\left[\begin{smallmatrix}
-t^1t^5&\frac{-t^1t^5}{-1-t^4+t^3t^5}&\frac{t^1(t^5)^2}{-1-t^4+t^3t^5}&1+t^4&t^5&0\\
-t^1&\frac{-t^1}{-1-t^4+t^3t^5}&\frac{t^1t^5}{-1-t^4+t^3t^5}&t^3&1&1\\
t^1(-1-t^4+t^3t^5)&t^1&-t^1t^5&0&0&2+t^4\\
0&0&0&0&0&0
\end{smallmatrix}\right]$.

The submatrix consisting of the first three columns of $A$ has rank 1,
so we can transform this matrix into a simpler matrix.

Recall that we assume that $t^1\ne 0$. When $t^5\ne 0$, then it turns
out that $d^\infty\sim d_1(\alpha:\beta)$, where
$\frac{\alpha\beta}{(\alpha+\beta)^2}=\frac{1-t^4-t^3t^5}{(2+t^4)^2}$.
Also, if $t^5=0$ and $t^3\ne0$, then the deformation is equivalent to
$d_1(1+t^4:1)$. In particular, if $t^4=0$, we see that there is a jump
deformation from our codifferential to the codifferential $d_1(1:1)$.
On the other hand, if $t^3=0$ and $t^4\ne 0$, we also deform to
$d_1(1+t^4:1)$. When $t^3=t^4=t^5=0$, there is a jump deformation of
$d^\infty$ to $d_1^\sharp$.

The picture for this element is more intriguing than for $d_3(1:-2)$.
In addition to the usual deformations along the family
$d_3(\alpha,\beta)$, jump deformation to $d_3(1:1:2)$, and
deformations along the big family, we see that $d_3(1:2)$ has a jump
deformation to the codifferential $d_1^\sharp$.  Because $d_1^\sharp$
itself has a jump deformation to $d_1(1:1)$, we get a jump deformation
to this element as well, and deformations along the family
$d_1(\lambda:\mu)$.  Thus we pick up far more deformations than we
would expect considering that the dimension of $H^2$ is only one more
than in the generic case. Again, the explanation for this ``impossibility''
has to do with the fact that the three dimensional tangent space to this element of the moduli
space does not accurately reflect the nature of the deformations, which are all
tangent to one of three planes in this space.  The true picture is captured by
the versal deformation.

\subsection{The codifferential $d_3(0:1)=\mathfrak r_2(\C)\oplus\C^2$}
The cohomology is given by
\begin{align*}
H^1=\langle&\psa12,\psa21-\psa31,\psa22-\psa32,\psa41,\psa42,\psa11\rangle\\
H^2=\langle&\psi_1=-\psa{12}2+\psa{13}2,\psi_2=\psa{14}1,\psi_3=\psa{14}2,\\
&\psi_4=\psa{24}1,\psi_5=\psa{24}2,\psi_6=\psa{12}1-\psa{13}1\rangle\\
H^3=\langle&\pha{124}1,\pha{124}3\rangle.
\end{align*}
Not all of the brackets of the nontrivial 2-cocycles vanish, so we
obtain some relations on the base of the versal deformation.  The
second order relations are $t^1t^2+t^3t^6=0$ and $t^1t^4+t^5t^6=0$.
The relations on the base are obtained by adding higher order terms to
these second order relations. We will omit them for brevity, but
instead will describe the solutions which may give actual
deformations.  There are 8 solutions, 4 of which not local. The local
solutions are
\begin{align*}
1)&\quad t^1=t^6=0\\
2)&\quad t^1=t^5=t^3=0\\
3)&\quad t^2=t^3=0, t^4=\tfrac{-t^5t^6}{t^1}\\
4)&\quad t^6=\tfrac{t^1t^2}{t^3}, t^4=\tfrac{t^2t^5}{t^3}.
\end{align*}

The first solution corresponds to the matrix
$
A=\left[\begin{smallmatrix}
0&0&0&t^2&t^4&0\\
0&0&0&t^3&t^5&1\\
0&0&0&0&0&1\\
0&0&0&0&0&0
\end{smallmatrix}\right]$.
The codifferentials $d^\infty$ associated to this matrix are easy to
analyze.  They usually lie in the big family, except for some special
cases when they are in the small family.  There is a jump
deformation to $d_3(0:0:1)$.

The second solution corresponds to
$
A=\left[\begin{smallmatrix}
t^6&-t^6&\frac{t^4t^6}{t^2-1}&t^2&t^4&0\\
0&0&0&0&0&1\\
0&0&0&0&0&1\\
0&0&0&0&0&0
\end{smallmatrix}\right]$.
If $t^6\ne0$, then either $t^2\ne1$ or $t^4=0$, and the deformation is
equivalent to $d_2^\sharp$. As a consequence, there is a jump
deformation to $d_2^\sharp$. When $t^6=0$, then if $t^4\ne0$ or
$t^2\ne1$, then the deformation is equivalent to $d_3(1:t_2:0)$.

In the third solution, let us first assume $t^6$ does not vanish. Then
the solution has matrix
$
A=\left[\begin{smallmatrix}
-(t^5-1)t^6&(t^5-1)t^6&\frac{t^5t^6}{t^1}&t^5&\frac{t^5}{t^1}&0\\
0&0&0&0&0&1\\
0&0&0&0&0&1\\
0&0&0&0&0&0
\end{smallmatrix}\right]$.
When $t^5\ne1$, then the deformation is equivalent to $d_2^\sharp$, 
so there is a jump deformation to $d_2^\sharp$.

Now let us assume that $t^6=0$ in the third solution. We can assume
$t^1\ne1$, since that corresponds to the first solution. The matrix
simplifies to
$
A=\left[\begin{smallmatrix}
0&0&0&0&0&0\\
t^1(t^5-1)&t^1&0&0&t^5&1\\
0&0&0&0&0&1\\
0&0&0&0&0&0
\end{smallmatrix}\right]$.
If $t^5\ne 1$, then the deformation is equivalent to $d_2^\sharp$ again.

Finally, let us consider the fourth solution, whose matrix is equivalent to
$
A=\left[\begin{smallmatrix}
\frac{t^1(t^2+t^5-1)}{t^2-1}&\frac{t^1t^2(t^2+t^5-1)}{t^3(t^2-1)}&t^1&t^2+
t^5&t^3&1\\
0&0&0&0&0&\frac{-t^2}{t^3}\\
0&0&0&0&0&1\\
0&0&0&0&0&0
\end{smallmatrix}\right]$.
There are two special cases that need to be considered, when $t^2=0$,
in which case, the restriction $t^3\ne0$ does not apply, and the case
when $t^5=0$. (The case $t^5=1-t^2$ is not local.) In these cases, the
restriction $t^2\ne1$ does not apply. Let us first address these
special cases.

When $t^2=0$, if $t^1\ne 0$ and $t^5\ne 1$ then we get $d_2^\sharp$.
On the other hand, if $t^2=0$ and $t^1=0$, then if $t^3=t^5=0$, we get
$d_3(0:1)$. 



From now on, we deal with the general case, assuming that $t^5\ne0$,
$t^2\ne0$ $t^5\ne 1-t^2$. Then $t^2\ne 1$ and $t^3\ne0$.

When $t^1\ne0$, the deformation is equivalent to $d_2^\sharp$;
otherwise it is equivalent to $d_3(1:t^2+t^5:0)$.

The subcases are a bit tricky, but the same codifferentials keep
showing up, so the final analysis of the deformations of this
codifferential is not difficult. We either obtain a  jump deformation to $d_3(1:0:0)$
or to $d_2^\sharp$, or we obtain a deformation along the big or small families.

\subsection{The codifferential $d_3(1:0)=\mathfrak r_{3,1}(\C)\oplus\C$}
The cohomology is given by
\begin{align*}
H^1=&\langle\pha11,\pha12,\pha22+\pha33,\pha43,\pha21+\pha33\rangle\\
H^2=&\langle-\psa{12}1+\psa{24}2+\psa{34}4,\psa{12}2+\psa{14}4,\psa{23}1,\psa{34}3,\psa{24}1,\psa{24}2,\psa{14}3\rangle\\
H^3=&\langle\pha{124}1,\pha{124}3,\pha{234}1\rangle.
\end{align*}
Some of the brackets of the nontrivial 2-cocycles do not vanish, and we have the second order relations
\begin{equation*}
t^1t^4+2t^3t^7+2t^2t^5=0,\quad  t^2t^4-t^2t^6+t^1t^7=0,\quad t^1t^5+t^3t^4+t^3t^6=0.
\end{equation*}
We omit the long expressions for the seven relations on the base of the
miniversal deformation, but remark that $1+t^4+t^6$ appears in the
denominator of two of them, so there may be an obstruction to the
extension of an infinitesimal deformation to a formal one.

The solution to the relations is quite complex; however, if we confine
ourselves to solutions which are local, then we can reduce the problem
to 9 relatively simple cases.
\begin{align*}
1)\quad t^1&=(t^4+t^6)\sqrt{\tfrac{t^2t^3}{t^4(1+t^6)}}, t^7=\tfrac{t^3(t^4+t^6)}{t^1},
t^5=\tfrac{-t^3(t^4+t^6)}{t^1}\\
2)\quad t^3&=\tfrac{(t^1)^2((t^6-t^4+2)^2+2(t^4-t^6))}{8t^2(t^4+t^6+2)},
  t^7=\tfrac{(t^4-t^6)(t^4-t^6-2)t^1}{8t^3},
  t^5= -\tfrac{t^3(t^4+t^6)}{t^1}\\
3)\quad t^1&=t^2=t^3=0\\
4)\quad t^1&=t^2=t^7=0,\quad t^6=-t^4\\
5)\quad t^1&=0, t^4=-t^6,t^7=\tfrac{t^6(1+t^6)}{t^5},
  t^2=\tfrac{-t^3t^6(1+t^6)}{(t^5)^2)}\\
6)\quad t^1&=t^3=t^5=0\\
7)\quad t^2&=t^4=t^7=0,\quad t^5=\tfrac{-t^3t^6}{t^1}\\
8)\quad t^4&=t^5=t^6=t^7=0\\
9)\quad t^3&=t^4=t^5=0,\quad t^7=\tfrac{t^2t^6}{t^1}.
\end{align*}
In the first solution, if $t^6=0$, or $4t^4\ne(t^6-t^4)^2$ then the
deformation is equivalent to $d_2^\sharp$; otherwise, we get
$d_1(1:0)$.

In the second solution   the differential is equivalent to
$d_1(\alpha:\beta)$, where
\begin{equation*}
(\alpha,\beta)=t^4+t^6+2\pm\sqrt{5(t^4)^2-12t^4-6t^4t^6+4t^6+5(t^6)^2+4}.
\end{equation*}
(It may be more revealing to recognize this element as $\mathfrak
g_8\left(\tfrac{4t^4-(t^6-t^4)^2}{(t^4+t^6+2)^2}\right)$). Note that
since $t^1$ is any nonzero number, this means that there is a jump
deformation from $d_3(1:0)$ to $d_1(1:0)$.

In the third solution, all deformations are either along the big
family or the family $d_3(\alpha:\beta)$. If $t^4=t^6=t^7=0$ and
$t^5\ne0$, then the deformation is equivalent to $d_3(1:1:0)$, so
there is a jump deformation from $d_3(1:0)$ to this element.

In the fourth solution, we get $d_1(1+t^6:-t^6)$, so that if $t^6=0$,
we see that there is a jump deformation to $d_1(1:0)$.

In the fifth solution, when $t^6=0$ then if $t^3=0$, we get
$d_3(1:1:0)$; otherwise we obtain $d_1(1:0)$. If $t^3=0$ and
$t^6\ne0$, then we get $d_3(1+\sqrt{t^6(t^6+1)}:1-\sqrt{t^6(t^6+1)})$
(assuming $t^6\ne-1$). When neither $t^6$ nor $t^3$ vanish, we get
a jump deformation to $d_2^\sharp$.

In the sixth solution, if $t^2=0$, this reduces to a previous case. If
$t^6=0$, then we get $d_2^\sharp$, a jump deformation. Otherwise, when
$t^7\ne0$ or $t^6\ne 1$, we get $d_1(1:t^6)$.

In the seventh solution, we always get $d_2^\sharp$.

In the eight solution, if $t^1=t^2=0$, this is a previous case. If
$t^2=0$, but $t^1\ne0$, or $t^1=0$ and $t^2\ne 0$, we get $d_2^\sharp$
unless $t^3=0$, in which case we get $d_1(1:0)$.  When neither $t^1$
nor $t^2$ vanish, then we get $d_2^\sharp$; unless $(t^3)^2=4t^1t^2$,
when we get $d_1(1:0)$.

In the ninth solution, we always get $d_2^\sharp$.

To summarize, we note that $d_3(1:0)$ jumps to $d_3(1:0:0)$ and $d_1(1:0)$
and it deforms along the the big
and small families as usual. It also jumps to $d_2^\sharp$.

\subsection{The codifferential $d_3^\star$}
 The cohomology is given by
\begin{align*}
H^1=&\langle\pha12,\pha13,\pha21,\pha22,\pha23,\pha31,\pha32,\pha33\rangle\\
H^2=&\langle\psa{14}2,\psa{14}3,\psa{34}3,\psa{24}3,\psa{14}1,\psa{34}2,
    \psa{34}1,\psa{24}1\rangle.
\end{align*}
The brackets of all 2-cocycles with each other vanish, so the
infinitesimal deformation is miniversal. The matrix of $d^\infty$ is
given by
$
A=\left[\begin{smallmatrix}
0&0&0&1+t^5&t^8&t^7\\
0&0&0&t^1&1&t^6\\
0&0&0&t^2&t^4&1+t^3\\
0&0&0&0&0&0
\end{smallmatrix}\right]$.
The deformations are easy to analyze, because they are given by the
equivalence classes of similar matrices of the $3\times 3$ submatrix
given by the parameters.  It is easy to see that there are jump
deformations to $d_3(1:1:1)$ and $d_3(1:1)$, as well as deformations
along the families these two codifferentials belong to. There are no
other possibilities for local deformations.

\subsection{The codifferential $d_2^\star=\mathfrak n_4(\C)$}
The cohomology is given by
\begin{align*}
H^1=&\langle\pha43,2\pha11+\pha22+\pha44,\pha31,\pha11+\pha22+\pha33\rangle\\
H^2=&\langle\psa{24}2,\psa{13}1+\psa{23}2,\psa{24}3,\psa{23}4,
\psa{14}3,\psa{12}1+\psa{13}2+\psa{23}3\rangle\\
H^3=&\langle\pha{124}2,\pha{124}3,\pha{123}1,\pha{123}4,
\pha{124}4-\pha{123}3\rangle.
\end{align*}
With such a large $H^3$, it would be too much to imagine that the brackets
of the cocycles vanish; in fact, there are 5 relations on the
base of the miniversal deformation. Since they are fairly simple, we
will give them:
\begin{align*}
&4t^2t^5+(t^1)^2t^6=0\\
&2t^5t^6-t^1t^2t^5+t^1t^3t^6=0\\
&t^1t^4+t^2t^6=0\\
&2t^2t^3t^4-t^1t^4t^6-t^1(t^2)^2t^6=0\\
&2t^4t^5+2t^1t^3t^4-(t^1)^2t^2t^6+t^1(t^6)^2=0.
\end{align*}
Note that the fourth relation has no second order term, and the fact
that the relations have no denominators means that the miniversal
deformation is constructed in a finite number of steps; in fact, since
the highest degree term in a relation is of degree 4, the fourth order
deformation is miniversal.  The solution to the relations can be
decomposed into 3 four dimensional subspaces and one more complex four
dimensional piece as follows.
\begin{align*}
1)\quad& t^4=t^5=t^6=0\\
2)\quad& t^2=t^4=t^6=0\\
3)\quad& t^1=t^2=t^5=0\\
4)\quad& t^6=\tfrac{-t^1t^4}{2t^2},\quad t^3=
  \tfrac{-(t^1)^2(t^4+(t^2)^2)}{4(t^2)^2},
\quad t^5=\tfrac{t^4(t^1)^3}{8(t^2)^2}.
\end{align*}
For the first solution, the matrix of $d^\infty$ is
$
A=\left[\begin{smallmatrix}
0&t^2&0&0&1&0\\
t^2t^3&0&t^2&0&t^1&1\\
0&t^2t^3&0&0&t^3&0\\
0&0&0&0&0&0
\end{smallmatrix}\right]$.  When $t^2=0$,
$d^\infty\sim d_3\left(0:\tfrac{t^1+\sqrt{(t^1)^2+4t^3}}{2}:
  \tfrac{t^1-\sqrt{(t^1)^2+4t^3}}{2}\right)$.
Assume $t^2\ne0$. Then, when $t^3=0$, if $t^1\ne0$, we have $d^\infty\sim
d_2^\sharp$, and when $t^1=0$, we get $d_1(1:0)$. When $t^3\ne0$, if
$(t^1)^2+4t^3=0$, then $d^\infty\sim d_1(1:0)$; otherwise it is
equivalent to $d_2^\sharp$. Thus we get jump deformations to
$d_2^\sharp$ and $d_1(1:0)$.

For the second solution, the matrix is given by
$
A=\left[\begin{smallmatrix}
0&0&0&0&1&0\\
0&0&0&0&t^1&1\\
0&0&0&t^5&t^3&0\\
0&0&0&0&0&0
\end{smallmatrix}\right]$.

In this case $d^\infty\sim d_3(\alpha:\beta:\eta)$ where
\begin{equation*}
\alpha+\beta+\eta=t^1q,\quad \alpha\beta\eta=t^5q^3,\quad \alpha\beta+
\alpha\eta+\beta\eta=t^3q^2,
\end{equation*}
where $q$ is an arbitrary nonzero parameter.  As a consequence, there
is a jump deformation from $d_2^*$ to every member of the family
$d_3(\alpha:\beta:\eta)$.

For the third solution, the matrix is given by
$
A=\left[\begin{smallmatrix}
t^6&0&0&0&1&0\\
0&t^6&0&0&0&1\\
t^6t^3&0&t^6&0&t^3&0\\
0&0&t^4&0&0&0
\end{smallmatrix}\right]$.

When $t^6=0$, if $t^3=0$, then we get $d_1(1:-1)$, while if $t^4=0$ we
get $d_3(1:-1:0)$, both jump deformations. When $t^6=0$ and neither
$t^3$ nor $t^4$ vanishes, then the deformation is equivalent to $d_3$,
another jump deformation.

Assume $t^6\ne0$. If $(t^6)^2\ne -t^3t^4$ then we get $d_3$; otherwise
we get $d_1(1:-1)$, both jump deformations.

For the fourth solution, the matrix is quite complicated, so we omit
it. When $t^1=0$, then if $t^2=0$, we get $d_1(1:-1)$,  and if
$t^4=0$, then we get $d_1(1:0)$, both jump deformations; otherwise, we
get a deformation along the family $d_1(\alpha:\beta)$.

When $t^1\ne0$ and $t^4=0$, then if $t^2=0$, we get the jump
deformation $d_3(1:1:0)$; otherwise we get a jump deformation to
$d_1(1:0)$.

When $t^1\ne0$ and $t^4\ne0$ and $t^4=-(t^2)^4$, then we get a jump
deformation to the element $d_1(1+\sqrt5,1-\sqrt5)$, which is just
$\mathfrak g_8(-1)$ on the Burde-Steinhoff list.

When none of the three conditions above hold, then the deformation is
equivalent to
$d_1(t^2+\sqrt{(t^2)^2-4t^4},t^2-\sqrt{(t^2)^2-4t^4})=\mathfrak
g_8\left(\tfrac{t^4}{(t^2)^2}\right)$. Since $t^1$ is an arbitrary
nonzero number, these deformations are also jump deformations.  Thus
there is a jump deformation to any element of the family
$d_1(\alpha:\beta)$.

To summarize, the deformations of $d_2^*$ are as follows. There are
jump deformations to every member of the big family
and everything they deform to, which means we get
jump deformations to the elements $d_2^\sharp$, $d_3$ and every element in the
family $d_1(\lambda:\mu)$.

\subsection{The codifferential $d_1=\mathfrak n_3(\C)\oplus\C$}
The cohomology is given by
\begin{align*}
H^1=\langle&\pha11+\pha22,\pha11+\pha44,\pha24,\pha31,\pha33,\pha42,\pha43\rangle\\
H^2=\langle&-\psa{23}3,\psa{23}4,\psa{14}1,\psa{14}2,\psa{14}3,
\psa{12}3,\psa{34}3,\\
&\psa{24}4,\psa{34}2\psa{12}4,\psa{13}1+\psa{23}2,
\psa{14}4-\psa{12}2,\psa{34}4-\psa{13}1\rangle\\
H^3=\langle&\pha{234}2,\pha{234}4,\pha{124}4,\psa{134}2,\psa{134}3,
\psa{123}2+\psa{134}4,\psa{123}3,\psa{123}4,\psa{124}2,\psa{124}3\rangle.
\end{align*}
There are 10 relations, none of which involve terms of higher order
than 3, so in fact, the second order deformation is already
miniversal.  We will not give the relations here explicitly. Because
the miniversal deformation is obtained in a finite number of steps,
the relations are polynomial, not rational, in the parameters.

If  all the parameters but $t^{13}$ and $t^{11}$ vanish, and
$t^{13}=-t^{11}$ then the relations are satisfied, and we have a jump
deformation to $d_1^\sharp$.

If all the parameters but $t^{12}$, $t^5$ and $t^6$ vanish, then if
$t^{12}\ne0$, the deformation is equivalent to $d_3$, so there is a
jump deformation to $d_3$.

If we assume that
\begin{align*}
t^8=&t^3=0,\quad t^1=\tfrac{-t^6t^7}{t^5},\quad t^{12}=\tfrac{-t^6t^9}{t^5},\\
t^2=&\tfrac{(t^6)^2t^9}{(t^5)^2},\quad t^4=\tfrac{-t^{12}t^5}{t^6},\quad t^{11}=\tfrac{-t^6t^9}{t^5},
\quad t^{10}=\tfrac{t^6t^{12}}{t^5},
\end{align*}
then if $q$ is a free parameter and
\begin{align*}
&\alpha+\beta+\eta=\tfrac{t^7q}{t^5}\\
&(\alpha+\beta)(\alpha+\eta)(\beta+\eta)=\tfrac{-t^9q^3}{(t^5)^2}\\
&(\alpha\beta+\alpha\eta+\beta\eta)=\tfrac{t^{12}q^2}{t^5t^6},
\end{align*}
we obtain a solution to the relations for which $d^\infty\sim
d_3(\alpha:\beta:\eta)$.
Whenever $t^5$ and $t^6$ don't vanish, there is a solution for every
$(\alpha:\beta:\eta)$, which means that there is a jump deformation
from $d_1$ to every element in the big family.

If all the parameters vanish except $t^3$, $t^4$ and $t^7$, then we
obtain a solution for which $d^\infty\sim d_3(\alpha:\beta)$, where
\begin{equation*}
\alpha=t^7q,\quad\alpha+\beta=t^3q,\quad\alpha\beta=-t^4q^2,
\end{equation*}
where again, $q$ is a nonzero free parameter, so we also have jump
deformations to every member of this family.

Similarly, if all the parameters but $t^5$, $t^8$, $t^{10}$ and $t^1$
vanish, and $t^8=t^1$, the relations are all satisfied, and if $q$ is
a nonzero parameter, then independently of the value of $t^5$ we have
$d^\infty\sim d_1(\alpha:\beta)$, where $q$ is a nonzero parameter and
\begin{align*}
\alpha+\beta=&-t^8q\\
\alpha\beta=&t^{10}q^2,
\end{align*}
so there is a jump deformation from $d_1$ to every member of the
family $d_1(\alpha:\beta)$.

If $t^5$ and $t^8$ do not vanish, $t^1=t^{12}=t^8$ and all the other
parameters vanish, then the relations are satisfied and  $d^\infty\sim
d_2^\sharp$, which gives another jump deformation.

If all the parameters except $t^5$ vanish, the relations are
satisfied, and we get a jump deformation to $d_2^*$.

Thus finally, we observe that $d_1$ has jump deformations to every
codifferential except $d_3^*$.  This pattern is completely analogous
to the three dimensional Lie algebra case, where the corresponding
element $d_1$ has jump deformations to every element except $d_2$,
which is exactly the analog of the element $d_3^*$, one dimension
lower.

\section{Description of the Moduli Space}

\begin{figure}

\includegraphics[width=4.9in]{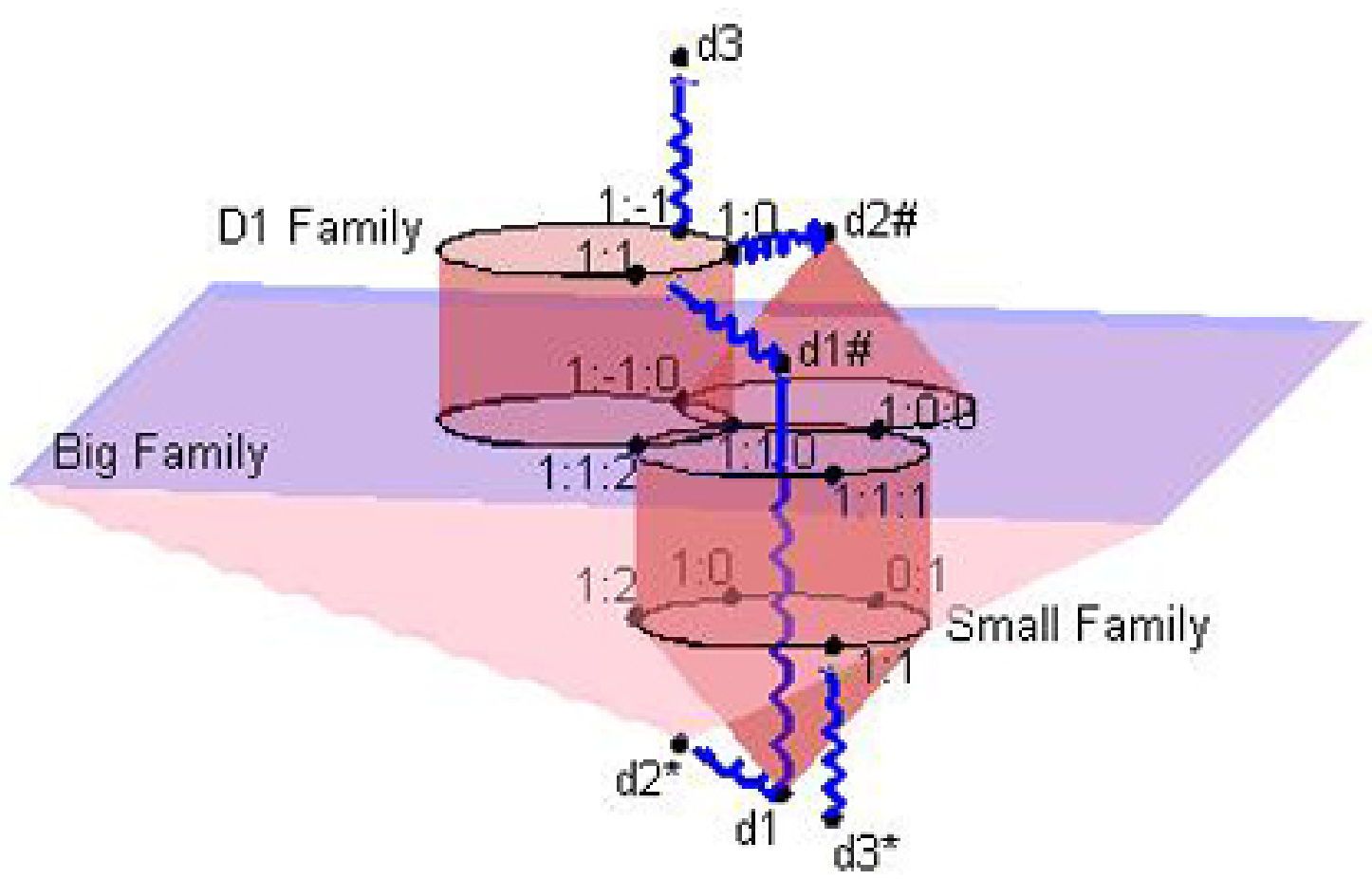}   

\caption{The Moduli Space of 4 dimensional Lie Algebras}
\label{modpic}
\end{figure}

In Figure (\ref{modpic}), we give a pictorial representation of the moduli space.
The big family $d_3(\lambda:\mu:\nu)$ is represented as a plane, although in reality it is
$\P^2/\Sigma_3$.
The families $d_1(\lambda:\mu)$, $d_3(\lambda:\mu)$ and the three
subfamilies $d_3(\lambda:\mu:0)$, $d_3(\lambda:\lambda:\mu)$ and $d_3(\lambda:\mu:\lambda+\mu)$
are represented by circles, mainly to reflect that the three subfamilies of the big family
 intersect in more than one point, because they each represent not
a single $\P^1$, but several copies of $\P^1$ which are identified under the action of the
symmetric group.

In the picture, jump deformations from special points are represented by curly arrows.
The jump deformations from the small family $d_3(\lambda:\mu)$ to $d_3(\lambda:\lambda:\mu)$
 and the jump deformations from
$d_3(\lambda:\mu:\lambda+\mu)$ to $d_1(\lambda:\mu)$ are represented by cylinders.
The jump deformations from the family $d_3(\lambda:\mu:0)$ to $d_2^\sharp$ and those
from $d_1$ to the small family are represented by cones.  Finally, the jump deformations
from $d_2^*$ to the big family are represented by an inverted pyramid shape.  All jump
deformations are either in an upwards or a horizontal direction.

The picture tries to capture the order of precedence of the deformations.  For example,
in the picture, you can trace a path of jump deformations from
$d_1$ to $d_3(1:0)$ to $d_3(1:1:0)$ to $d_1(1:0)$ to $d_2^\sharp$.

\section{Classifying a Particular Lie Algebra}

In \cite{Aga}, it was shown that a four dimensional Lie algebra
can be classified by computing certain invariants of the Lie algebra. Instead,
our approach to classifying a Lie algebra, which we will
outline here, used linear algebra.

Suppose that a codifferential $d$ representing a Lie algebra structure has matrix $A$.
Since the rank of the matrix is at most 3, it is easy to compute a new basis for which
the matrix has the form $A=\left[\begin{smallmatrix}A'&\delta\\0&0\end{smallmatrix}\right]$,
where $A'$ is a $3\times 3$ matrix representing a 3 dimensional Lie algebra, and $\delta$
is a $3\times 3$ matrix representing a derivation of this Lie algebra.

Next, consider the submatrix $A'$. If it has rank 3, then the codifferential is equivalent
to $d_3$.  Otherwise, we find a new basis in which the submatrix $A'$ has been reduced to
one with exactly as many rows as its rank. In fact, by using the classification methods
for three dimensional Lie algebras, one can reduce the matrix $A'$ to one of the standard
forms.

At this point, the matrix $\delta$ representing the derivation on the three dimensional Lie
algebra may not represent an outer derivation. However, by replacing the vector $e_4$ with
a vector of the form $e_4'=ae_1+be_2+ce_3+e_4$, one can replace the $\delta$ with one
representing an outer derivation.

Once this has been accomplished, the classification scheme presented in this paper for
determining the point in the moduli space corresponding to an extension of a three dimensional
Lie algebra by an outer derivation can be applied.  The precise identification scheme depends
on which point in the moduli space of three dimensional Lie algebras occurs.

When computing versal deformations of the four dimensional Lie algebras, in most cases,
we could identify the appropriate element by following a more simple scheme of solving for
a matrix $G$ such that the matrix $GA'= A Q$, where $Q$ is the matrix representing the
linear transformation $g$ corresponding to $G$ extended to $\E^2V\ra\E^2V$ and
$A'$ is a matrix representing one of the nine types of elements in the moduli space.

However, because our matrices involved many parameters, it was sometimes
too difficult for the
computer to solve for the values of the parameters for which the $A'$ and $A$ matrices
are equivalent.  In those cases, we followed the more complicated scheme outlined above.
In practice, we found that it was only necessary to follow the steps partially, because
after transforming the matrix to eliminate some of the rows, we then were able to apply
the simple scheme, and obtain a solution.

\section{Conclusions}
The computation of the equivalence classes of non-isomorphic Lie
algebra structures in a vector space $V$ determines the elements of
the moduli space of Lie algebra structures on $V$, but is only the
first step in the classification of these structures. When
classifying the algebras, there are different ways of dividing up
the structures according to families; therefore, it is desirable to
have a rationale for the division.  In this paper, we have shown
that there is a natural way to divide up the moduli space into
families, using cohomology as a guide to the division, and versal
deformations as a tool to refine the analysis.

The four dimensional Lie algebras can be decomposed into families,
each of which is naturally an orbifold.  If one takes into account
the information about jump deformations,  the division we have given
is uniquely determined. The elements of the family which contain a
Lie algebra structure $d$ are precisely those Lie algebras which can
be obtained as smooth deformations of $d$, but which are not smooth
deformations of any Lie algebra structure $d'$ which is a jump
deformation of $d$. This rule allows us to distinguish between the
algebra $d_3^*$ and $d_3(1:1)$, for example. Even though $d_3^*$ has
smooth deformations to the family $d_3(\lambda:\mu)$, it also has a
jump deformation to $d_3(1:1)$, which has smooth deformations to the
same family.  Thus $d_3(1:1)$, which has no jump deformations to any
element which has smooth deformations to the family, is the element
which belongs to the family.

According to this system,  there is one two-parameter family, two
one-parameter families, and six singleton elements, giving rise to a
two-dimensional orbifold, two one-dimensional orbifolds, and six
one-dimensional orbifolds.  The jump deformations provide maps
between the families which either are smooth maps of orbifolds (or
suborbifolds as in the case of the map
$d_3(\lambda:\mu:\lambda+\mu)\ra d_1(\lambda:\mu)$), or, in the case
of some of the singletons, identify the element with a whole family.

The cohomology of a Lie algebra determines the tangent space to the
Lie algebra, but the tangent space does not contain enough
information to give a good local description of the moduli space.
The relations on the base of the versal deformation determine the
manner in which the moduli space contacts the tangent space. In one
example here, the tangent space was two dimensional, but
deformations were only along two curves.  In another case, the
tangent space was three dimensional, but the deformations were
confined to three planes. It is clear that the cohomology is not
sufficient to get an accurate picture of the moduli space. 
Versal deformations provide important detail that characterizes the
moduli space completely.

\section{Acknowledgements}
The authors would like to thank E. Vinberg for helpful discussions and the
Max-Planck-Institut f\"ur Mathematik, Bonn for hosting both authors while they were
finishing this paper.
\bibliographystyle{amsplain}

\begin{thebibliography}{10}

\bibitem{Aga}
Y.~Agaoka, \emph{An algorithm to determine the isomorphism classes of
  4-dimensional complex {L}ie algebras}, Linear Algebra and its Applications
  \textbf{345} (2002), 85--118.

\bibitem{bs}
D.~Burde and C.~Steinhoff, \emph{Classification of orbit closures of
  4-dimensional complex {L}ie algebras}, Journal of Algebra \textbf{214}
  (1999), 729--739.

\bibitem{fi1}
A.~Fialowski, \emph{Deformations of {L}ie algebras}, Mathematics of the
  USSR-Sbornik \textbf{55} (1986), no.~2, 467--473.

\bibitem{fi2}
A.~Fialowski, \emph{An example of formal deformations of Lie algebras},
Proc. NATO Conf. on Deformation Theory of Algebras and Appl., Kluwer
1988, 375--401.

\bibitem{ff2}
A.~Fialowski and D.~Fuchs, \emph{Construction of miniversal deformations of
  {L}ie algebras}, Journal of Functional Analysis (1999), no.~161(1), 76--110.

\bibitem{fp1}
A.~Fialowski and M.~Penkava, \emph{Deformation theory of infinity algebras},
  Journal of Algebra \textbf{255} (2002), no.~1, 59--88, math.RT/0101097.

\bibitem{fp3}
\bysame, \emph{Versal deformations of three dimensional {L}ie algebras as
  \linf\ algebras}, Communications in Contemporary Mathematics \textbf{7}
  (2005), no.~2, 145--165, math.RT/0303346.

\bibitem{gers}
M.~Gerstenhaber, \emph{On the deformations of rings and algebras {I}--{IV}},
  Annals of Mathematics \textbf{79} (1964), 59--103; \ {II}, Annals of
  Mathematics \textbf{84} (1966), 1--19; \ {III}, Annals of Mathematics
  \textbf{88} (1968), 1--34; \ {IV}, Annals of Mathematics \textbf{99} (1974),
  257--276.

\bibitem{KN}
A.A.~Kirillov and Y.A.~Neretin, \emph{The variety $A_n$ of n-dimensional Lie
  algebra structures}, Amer. Math. Soc. Transl. \textbf{137} (1987), no.~2,
  21--30.

\bibitem{KT}
B.~Komrakov and A.~Tchourioumov, \emph{Small dimensional and linear {L}ie
  algebras}, International Sophus Lie Centre Press, 2000.

\bibitem{nr}
A.~Nijenhuis and R.~Richardson, \emph{Deformations of Lie algebra structures},
  Jour. Math. Mech. \textbf{17} (1967), 89--105.

\bibitem{ov}
A.L.~Onishik and and E.B.~Vinberg (ed.) \emph{Lie Groups and Lie Algebras
III.(Structure of Lie groups and Lie algebras)}, Encyclopaedia of Math. Sci.,
1994. 

\bibitem{pbnl}
R.O. Popovich, V.M. Boyko, M.O. Nesterenko, and M.W. Lutfullin,
  \emph{Realization of real low-dimensional {L}ie algebras}, Journal of Phys. A
  Math Gen. \textbf{36} (2003), 7337--7360.

\bibitem{Rh}
M.~Rhomdani, \emph{Classification of real and complex nilpotent {L}ie algebras
  of dimension 7}, Linear and Multilinear Algebra \textbf{24} (1989), 167--189.

\bibitem{Tu}
P.~Turkowski, \emph{Literature on the structure of low-dimensional
  nonsemisimple {L}ie algebras and its applications to cosmology}, Acta
  Cosmologica \textbf{20} (1994), 147--153.

\end{thebibliography}

\providecommand{\bysame}{\leavevmode\hbox to3em{\hrulefill}\thinspace}
\providecommand{\MR}{\relax\ifhmode\unskip\space\fi MR }
\providecommand{\MRhref}[2]{%
  \href{http://www.ams.org/mathscinet-getitem?mr=#1}{#2}
}
\providecommand{\href}[2]{#2}

\end{document}